%% file: jamesnew.tex
%%%%%%%%%%%%%%%%%%%%%%%%%% jamesnew.tex %%%%%%%%%%%% 
%
%  This is  FRS paper II (James bundles) edited for publication Jan 2003 
%
\input gtmacros
\input pictex
\input newinsert

\let\ninepoint\small
\hoffset 0.5truein     %  To get the
\voffset 1truein       %  print centrepage

%%%% FONTS :
\font\spec=cmtex10  %%% for \d (better \partial !)

%%%% SYMBOLS abbreviations etc.:
\def\bar{\overline}
\def\text{\hbox}
\def\a{\longrightarrow}
\def\inv{^{-1}}
\def\ep{\epsilon}
\def\d{\hbox{\spec \char'017\kern 0.05em}} %%%% better \partial %%%%

\def\ld{\ldots}
\def\la{\lambda}
\def\ga{\gamma}
\def\de{\delta}

\def\mo{\hbox{\hskip 1pt\vrule height 8.5pt depth 2pt width 0.7pt\hskip 1pt}}
\let\\=\cr

\def\wreath{\hbox{\tenmsa\lower0.4ex\hbox{\char'161}\kern-0.13em\raise0.2ex\hbox{\char'170}}}

\def\N{{\Bbb N}} 
 
\def\C{{\cal C}}
\def\D{{\cal D}}

\def\N{{\cal N}}

\def\Z{{\Bbb Z}}

\def\re{{\Bbb R}}

\def\ha{{\scriptstyle{1\over2}}}
\def\cupprod{{\scriptscriptstyle\,\cup\,}}

\def\tsq{\sqr57}
\def\ssq{\sqr34}
\def\nsq{\sqr45}

\mathchardef\square="0\hexa03
\mathchardef\celt="0\hexa01
\mathchardef\twid="1218
\mathchardef\diam="027D       %% diamond
\mathchardef\clock="0\hexa08
\mathchardef\aclock="0\hexa09
\mathchardef\rar="0225        %% arrow going up right
\mathchardef\lar="022D        %% arrow going up left
\mathchardef\ldar="022E       %% down left
\mathchardef\rdar="0226       %% down right

\def\r#1{\raise6pt\hbox{#1}}

%%% for labelled remarks

%%%  macros for importing figures from xfig  %%%
\long\def\mbox#1{$$\vbox{#1}$$}
\newdimen\unitlength
\unitlength= 1.000cm
\def\makebox(0,0)[l]#1{\hbox{#1}}
\def\symbol#1{\char#1}
\font\tencirc=lcircle10
\def\SetFigFont#1#2#3#4{{\small$#4$}}
\font\thinlinefont=cmr7
\def\co{\kern.1ex\colon\thinspace}

%%%%%%%%%%%%          Main text starts here        %%%%%%%%%
%
%                       References first
%
\reflist

\refkey\Ant {\bf R Antolini}, {\it PhD Thesis} Warwick (1996)

\refkey\AntWie {\bf R Antolini}, {\bf B Wiest}, {\it The homotopy
type of the singular cubical set of a topological space}, Warwick 
University Preprint (1997)

\refkey\Baues {\bf Baues}, {\it Commutator Calculus and Groups of Homotopy Classes}, London Math. Soc. Lecture Note Series, no. 50 C.U.P. (1981)

\refkey\Bruck {\bf R\,H Bruck}, {\it A survey of binary systems},
Springer, Berlin, 1958 

%\refkey\BLin 
%{\bf J\,S Birman}, {\bf X-S Lin}, {\it Knot polynomials and Vassiliev's
%invariants}, Invent. Math. 111 (1993) 225--270.

\refkey\BarNat {\bf D Bar-Natan}, {\it On the Vassiliev knot invariants}, 
Topology, 34 (1995) 423--472

\refkey\BRS 
{\bf S Buoncristiano}, {\bf C Rourke}, {\bf B Sanderson}, {\it A geometric 
approach to homology theory}, London Math. Soc. Lecture Note Series, no. 
18 C.U.P. (1976)

\refkey\Dre {\bf W Dreckmann},
{\it Distributivgesetze in der Homotopietheorie}, PhD Thesis, Bonn (1992)

\refkey\FeRo 
{\bf R Fenn}, {\bf C Rourke}, {\it Racks and links in codimension two},
Journal of Knot theory and its Ramifications, 1 (1992) 343--406

\refkey\Trunks {\bf R Fenn}, {\bf C Rourke}, {\bf B Sanderson}, {\it Trunks and
classifying spaces}, Applied Categorical Structures, 3 (1995)
321--356

\refkey\James {\bf R Fenn}, {\bf C Rourke}, {\bf B Sanderson}, 
{\it James bundles and applications}, preprint (1996) {\tt
http://www.maths.warwick.ac.uk/\char'176cpr/ftp/james.ps}

\refkey\Rackspace {\bf R Fenn}, {\bf C Rourke}, {\bf B Sanderson}, 
{\it The rack space}, (to appear)

\refkey\Class {\bf R Fenn}, {\bf C Rourke}, {\bf B Sanderson}, {\it A 
classification of classical links}, (to appear)

\refkey\Joyce {\bf D Joyce}, {\it A classifying invariant of knots; the knot quandle},
J. Pure Appl. Alg. 23 (1982) 37--65

%%%%%% for reduced products
\refkey\JRed {\bf I James}, {\it Reduced product spaces}, Ann. Math. 62
(1955) 170--197     

%%%%%% for the p-th James-Hopf invariants
\refkey\JHopf {\bf I James}, {\it On the suspension triad},  Ann. Math. 63
(1956) 191--247 

\refkey\Kont {\bf M Kontsevich}, {\it Vassiliev's knot invariants},
 I.M.Gelfand Seminar, 16, Part 2, 137--150, Adv. Soviet Math., 
Amer. Math. Soc., Providence RI, 1993

\refkey\KS {\bf U Koschorke}, {\bf B Sanderson}, {\it Self-intersections and
higher Hopf invariants}, Topology, 17 (1978) 283--290

%%%%% Immersions can be assumed self transverse
\refkey\LS {\bf R Lashof}, {\bf S Smale}, {\it Self-intersections of immersed
manifolds}, J. Math. Mech. 8 (1959) 143--157

\refkey\Mah {\bf M Mahowald}, {\it Ring Spectra which are Thom complexes},
 Duke Math. Journal, 46, no. 3, (1979) 549--559

\refkey\Priddy {\bf S Priddy}, {\it K(Z/2) as a Thom spectrum}, 
Proc. Amer. Math. Soc. (2) 70 (1978) 207--208

\refkey\Puppe {\bf D Puppe}, {\it Some well known weak homotopy equivalences
are genuine homotopy equivalences}, Symposia Mathematica, Vol.V
(INDAM, Rome, 1969/70), 363--374, Academic Press, London, 1971

\refkey\Delt {\bf C Rourke}, {\bf B Sanderson}, {\it $\Delta $-sets, I:
homotopy theory}, Quart. J. Math. (Oxford) (2) 22 (1971) 321--338\hfil\break
{\it $\Delta$-sets, II: block bundles and block fibrations}, 
Quart. J. Math. (Oxford) (2) 22 (1971) 465--485 

\refkey\SMah {\bf B Sanderson}, {\it The geometry of Mahowald orientations},
Algebraic Topology, Aarhus. Springer Lecture Notes in Mathematics,
no. 763 (1978) 152--174

\refkey\Vass {\bf V\,A Vassiliev}, {\it Cohomology of knot spaces}, Theory of
Singularities and its Applications, 23--69, 
Amer. Math. Soc., Providence, RI, 1990.

\refkey\Wiest {\bf B Wiest}, {\it PhD Thesis} Warwick (1997)

\endreflist

\title{James bundles}
\author{Roger Fenn}
\address{Department of Mathematics, University of Sussex\\
Falmer, Brighton, BN1 9QH, UK\\\bigskip
{\sc Colin Rourke\\Brian Sanderson}\\\medskip
Mathematics Institute, University of Warwick\\
Coventry, CV4 7AL, UK}
\email{R.A.Fenn@sussex.ac.uk, cpr {\rm and} bjs@maths.warwick.ac.uk}

\abstract

We study cubical sets without degeneracies, which we call
$\nsq$--sets. These sets arise naturally in a number of settings and
they have a beautiful intrinsic geometry; in particular a $\nsq$--set
$C$ has an infinite family of associated $\nsq$--sets $J^i(C)$,
$i=1,2,\ldots$, which we call James complexes.  There are mock bundle
projections $p_i\co|J^i(C)|\to|C|$ (which we call James bundles)
defining classes in unstable cohomotopy which generalise the classical
James--Hopf invariants of $\Omega(S^2)$.  The algebra of these classes
mimics the algebra of the cohomotopy of $\Omega(S^2)$ and the
reduction to cohomology defines a sequence of natural characteristic
classes for a $\nsq$--set.  An associated map to $BO$ leads to a
generalised cohomology theory with geometric interpretation similar to
that for Mahowald orientation [\Mah, \SMah].

\endabstract



\maketitle

This paper is concerned with $\sq$--sets, pronounced ``square sets'',
which are cubical sets without degeneracies.  These arise naturally in
a number of settings and can be used to classify links of codimension
2.  Here we shall describe some of the instrinsic geometry of a
$\sq$--set and in particular we shall define an infinite family of
associated $\sq$--sets, which we call {\it James complexes}.  These
define a sequence of mock bundles with base the given $\sq$--set, the
{\it James bundles}, which have strong connection with classical
James--Hopf invariants.  They define natural characteristic classes
for a $\sq$--set and lead to associated generalised cohomology
theories with geometric interpretation similar to that for Mahowald
orientation [\Mah, \SMah].

In an earlier paper [\Trunks] we explored the abstract connections
between $\sq$--sets and categories through the concept of a {\it
trunk} which can be regarded both as a generalisation of the (cubical)
nerve of a category and as a $\sq$--set.  We defined a classifying
space ($\sq$--set) for a rack (the {\it rack space}) and used the
trunk formalism to establish connections with classical classifying
spaces (for categories, groups and crossed complexes).  In
[\Rackspace] we shall apply the results of [\Trunks] and this paper
by showing that the rack space is a classifying space for smooth links
(of codimension 2) with the first James bundle as classifying bundle.
Further in [\Class] we shall give a classification of classical links
using racks and the canonical class in $\pi_2$ of the rack space
determined by the link diagram.

This paper formed part of our January 1996 preprint ``James bundles
and applications'' [\James], which also contains a preliminary version
of [\Rackspace].

Further results on $\sq$--sets are to be found in [\Ant] and
[\AntWie].  In [\Ant] it is shown that the natural map $C\to T$ of a
$\sq$--set to the trivial $\sq$--set induces a homotopy equivalence
$\mo C\mo\simeq\Omega(S^2)\times R$ in the case when $C$ admits
degeneracies and $R$ is the realisation using degeneracies.  In
particular, the singular $\sq$--set ${\rm S}(X)$ of a space $X$ has
the weak homotopy type $\Omega(S^2)\times X$.  This result is given an
alternative proof in [\AntWie] in terms of the natural interpretation
of maps into $C$ as labelled diagrams, which comes from using
transversality, see [\Rackspace].  However it is not the case that all
naturally occuring infinite $\sq$--sets split a copy of $\Omega(S^2)$
up to homotopy type.  For example rack spaces do not in general split
in this way (see the computations of homotopy type of rack spaces in
[\Rackspace, \Wiest]).  Moreover a general $\sq$--set can have an
arbitrary weak homotopy type (corollary 1.3).

This paper is organised as follows.  In section 1 we define
$\sq$--sets, give examples and explain the connection with
$\Delta$--sets.  In section 2 we define the James complexes and James
bundles of a $\sq$--set and in section 3 we give the connection with
the classical James--Hopf invariants.  In section 4 we consider the
classes in stable cohomotopy and in cohomology determined by the James
bundles and show how these classes are pulled back from
Steifel--Whitney and Wu classes of $BO$.  In section 4 we use the map
to $BO$ to define the concept of a $\sq$--set orientation for a
manifold.  This leads to a generalised (co)-homology theory which has
a geometric interpretation close to that for Mahowald orientation
given in [\SMah].

\section{Basic definitions and examples}

A {\sl $\sq$--set} $C=\{C_n;\d^\epsilon_i\}$ consists of a
collection of sets $C_n$,
one for each natural number $n\geq0$, and functions
$\d_i^\epsilon\co C_n\to C_{n-1},\ 1\leq i\leq n,\ \epsilon\in
\{0,1\}$, called {\sl first order face maps}, 
satisfying the following {\sl face relations}: 
$$\d^\eta_{j- 1}\d^\epsilon_i=\d^\epsilon_i\d^\eta_j, \quad 1\leq i<
j\leq n,\quad\ \epsilon,\eta\in\{0,1\}.$$

There is a  notion of a {\sl $\sq$--map} $f\co C\to D$
between $\sq$--sets, namely a collection of maps $f_n\co C_n\to D_n$
for $n=0,1,\ldots$, which commute
with the face maps, ie $f_{n-1}\d^\epsilon_i=\d^\epsilon_if_n$.

Note that a $\sq$--set differs from a 
{\it cubical set}, which has a similar definition but also has
degeneracy maps;  for more detail on the connections between
cubical and $\sq$--sets, see [\Trunks].

There is also the notion of a {\sl $\sq$--space}.  This comprises
a collection $\{X_n\}$ of topological spaces and first order face maps (continuous
maps) $\d_i^\epsilon\co X_n\to X_{n-1},\ 1\leq i\leq n,\ \epsilon\in
\{0,1\}$ satisfying the same face relations as above.  

We shall give presently alternative definitions for $\sq$--sets
and maps based on cubes.  First we need to establish some terminology.

\sh{The category $\tsq$}

The {\sl $n$--cube} $I^n$ is the subset $[0,1]^n$ of $\re^n$.   

A {\sl $p$--face} of $I^n$ is a subset defined by choosing $n-p$
coordinates and setting some of these equal to 0 and the rest to 1.  
In particular there are $2n$
faces of dimension $n-1$ determined by setting $x_i=\epsilon$ where
$i\in\{1,2,\ldots,n\}$ and $\ep\in \{0,1\}$.

A 0--face is called a {\sl vertex} and corresponds to a point
of the form
$(\ep_1,\ep_2,\ldots,\ep_n)$ where $\ep_i=0$ or $1$ and
$i\in\{1,2,\ldots,n\}$.  The 1--faces are called {\sl edges} and the
2--faces are called {\sl squares}.

Let $p\leq n$ and let $J$ be a $p$--face of $I^n$.  Then there is
a canonical {\sl face map} $\lambda\co I^p\to I^n$, with $\lambda(I^p)=J$,
given in coordinate
form by preserving the order of the coordinates $(x_1,\ldots,x_p)$
and inserting $n-p$ constant coordinates which are either 0 or 1.  
If $\la$ inserts only $0$'s (resp.\ only $1$'s) we call it a {\sl front}
(resp.\
{\sl back}) face map.
Notice that any face map has a unique {\sl front--back} decomposition as $\la\mu$, say,
where $\la$ is a front face map and $\mu$ is a back face map. There is also a unique 
{\sl  back--front} decomposition.
There are $2n$ face maps defined by the $(n-1)$--faces which are denoted
$\delta_i^\ep\co I^{n-1}\to I^n$, and given by:
$$\delta_i^\ep(x_1,x_2,\ldots,x_{n-1})=
(x_1,\ldots,x_{i-1},\ep,x_{i},\ldots,x_{n-1}),\qquad\epsilon\in \{0,1\}.$$

The following relations hold:
$$
\de_i^\ep\de_{j-1}^\omega
=\de_j^\omega\de_i^\ep,\quad 1\leq i<j\leq n,\quad
\ep,\omega\in \{0,1\}.\leqno{\bf \label}$$

\rk{Definition} The category $\sq$ is the category whose objects are the $n$--cubes $I^n$ for $n=0,1,\ldots$ and whose morphisms are the face maps.
\ppar

\sh{$\tsq$--Sets and their Realisations}

We now give the second equivalent definitions of $\sq$--sets and maps.

A  {\sl $\sq$--set} is a functor 
$C\co\sq^{op}\to Sets$ where $\sq^{op}$ is the opposite category 
of $\sq$ and $Sets$ denotes the category of sets.

A {\sl $\sq$--map} between $\sq$--sets is a natural transformation.

We write $C_n$ for $C(I^n)$, $\lambda^*$ for $C(\lambda)$ and we
write $\d_i^\ep$ for $C(\delta_i^\ep)=(\delta_i^\ep)^*$.

The equivalence of the two definitions of $\sq$--sets and maps
follows from the formul\ae\ for the composition of face maps given
in equation 1.1 above.

The {\sl realisation\/} $\mo C\mo$ of a $\sq$--set $C$ 
is given by making the identifications
$(\lambda^*x,t)\sim (x,\lambda t)$ in the disjoint union
$\coprod_{n\ge0} C_n\times I^n$. 
%We write $\mo c\mo$
%for the image of $\{c\}\times I^n$ in $\mo C\mo$.

We shall call 0--cells (resp.\ 1--cells, 2--cells) of $\mo C\mo$
vertices (resp.\ edges, squares) and this is consistent with the
previous use for faces of $I^n$, since $I^n$ determines a $\sq$--set
with cells corresponding to faces, whose realisation can be
identified in a natural way with $I^n$.

Notice that $\mo C\mo$ is a CW complex with one $n$--cell for each
element of $C_n$ and that each $n$--cell has a canonical
characteristic map from the $n$--cube.

However, not every CW complex with cubical characteristic maps
comes from a $\sq$--set --- even if the cells are glued by isometries
of faces.  In $\mo C\mo$, where $C$ is a $\sq$--set, cells are
glued by face {\it maps}, in other words by canonical isometries
of faces.  

There is a similar alternative definition of a $\sq$--space
namely a functor $X\co \sq^{op}\to Top$  (where $Top$ denotes the category
of topological spaces and continuous maps) and the realisation $\mo X\mo$
given by the same formula as above.  Note that a $\sq$--{\sl set} is effectively
the same as a $\sq$--{\sl space} in which all spaces have the discrete topology.

\sh{Connection with $\Delta$--sets}

The concept of a $\Delta$--set is similar to the concept of a
$\sq$--set but is based on simplexes rather than cubes.  A
$\Delta$--set can have an arbitrary (weak) homotopy type,
see [\Delt] where basic material on $\Delta$--sets can be
found.  The definition of a $\Delta$--set can be obtained
from the first definition of a $\sq$--set (at the start
of the section) by
forgetting the $\epsilon$'s and $\eta$'s and allowing $i=0$.  

A $\sq$--set can be subdivided to form a $\Delta$--set by coning
from centres.  More precisely suppose that
$C$ is a $\sq$--set. Define the $\Delta$--subdivision, $Sd_\Delta(C)$, 
as follows. A $k$--simplex of $Sd_\Delta(C)$ is a $(k+1)$--tuple
$(c,\la_1,\ldots,\la_k)$ where $c\in C_n$, for some $n$,
$\la_i\co I^{n_{i-1}}\to I^{n_i}$ is a face map, $1\le i\le k$,
 and $n_0<n_1<\dots <n_k=n$.
The face maps $\d_i$, for $0\le i\le k$, are given by
$$\d_i(c,\la_1,\ldots,\la_k)=\cases{(c,\la_2,\ldots,\la_k)
&for $i=0$\cr
(c,\la_1,\ldots,\la_{i-1},\la_{i+1}\la_i,\la_{i+2},\ldots,\la_k)
&for $0<i<k$\cr
(\la_k^*c,\la_1,\ldots,\la_{k-1})&for $i=k$\cr}$$
Conversely given a $\Delta$--set $X$
we can define a $\sq$--subdivision, $Sd_\ssq(X)$ of $X$ by projecting 
simplexes from the origin.   More precisely, let 
$\Delta^n=\{x\in \re^{n+1}\mid x_i\ge0,\sum_i x_i=1\}$ denote the 
standard  n--simplex.   Let $\rho_n\co
\Delta^n\to\d I^{n+1}$ be the radial projection from $0\in \re^{n+1}$.
Then $\rho_n(\Delta^n)$ is the union of the back $n$--faces of $I^{n+1}$.
A picture of the image of $\Delta^2$
under this projection (seen from the origin) is
given in figure \figkey\SubDiv, and it can be seen that $\Delta^2$ has been subdivided
into three squares.  A similar projection subdivides
an $n$--simplex into $n+1$ $n$--cubes and this process defines the
required $\sq$--subdivision.

\fig{\SubDiv}
\beginpicture\ninepoint
\setcoordinatesystem units <.5truein,.5truein> point at 0 0
\setplotarea x from -1 to 1, y from 0 to 1.732
\setlinear
\plot 0 0 0 0.577 0.5 0.866 /
\plot 0 0.577 -0.5 0.866 /
\arrow <2.5pt> [0.3, 1] from -0.6 0 to -0.5 0
\arrow <2.5pt> [0.3, 1] from 0.6 0 to 0.5 0
\arrow <2.5pt> [0.3, 1] from 0 0.189 to 0 0.289
\arrow <2.5pt> [0.3, 1] from -1 0 to -0.75 0.433
\arrow <2.5pt> [0.3, 1] from 1 0 to 0.75 0.433
\arrow <2.5pt> [0.3, 1] from -0.5 0.866 to -0.25 0.7215
\arrow <2.5pt> [0.3, 1] from 0.5 0.866 to 0.25 0.7215
\arrow <2.5pt> [0.3, 1] from 0 1.732 to -0.25 1.299
\arrow <2.5pt> [0.3, 1] from 0 1.732 to 0.25 1.299
\setplotsymbol ({\tenrm .})
\setlinear
\plot -1 0 1 0 0 1.732 -1 0 /
\endpicture
\endfig

For detailed formul\ae, 
let  $B(k,n)$ denote the set of back face maps $\la\co I^k\to I^n$.
 If $\la\co I^k\to I^n$ is a front face map let
$r(\la)\co \Delta^{k-1}\to\Delta^{n-1}$ denote its
 restriction.  Define the $\sq$--subdivision, $Sd_{\ssq}(X)$ as follows. Set
$$Sd_{\ssq}(X)_k=\coprod_{n>k} X_{n-1}\times B(k,n).$$
Supose given $x\in X_{n-1},\>\la\in B(k,n)$ and
 a face map $\mu\co I^s\to I^k$. Let $\la'\mu'$ be the
 front--back decomposition of $\la\mu$.  Define
$\mu^*(x,\la)=(r(\la')^*x,\mu')$.

\proc{Proposition}\items
\item{\rm(1)}For any $\Delta$--set $X$ the spaces $\mo X\mo$
and $\mo Sd_{\ssq}(X)\mo$ are homeomorphic,
\item{\rm(2)}For any $\sq$--set $C$ the spaces $\mo C\mo$
and $\mo Sd_{\Delta}(C)\mo$ are homeomorphic.\enditems\endproc
\prf The homeomorphism $\mo Sd_{\ssq}(X)\mo\to \mo X\mo$ is induced
by 
$$(x,\la,t)\in X_{n-1}\times B(k,n)\times I^k \longmapsto
(x,\rho^{-1}(\la(t)))\in X_{n-1}\times\Delta_{n-1}.$$
Let $\hat c_m=(\ha,\ldots,\ha)\in I^m$.
The homeomorphism $\mo Sd_{\Delta}(C)\mo\to \mo C\mo$ is induced
by 
$$\eqalign{
((c,\la_1,\ldots,\la_k),{\bf t})\in Sd_\Delta(C)_k\times\Delta_k&
\longmapsto\cr
(c,t_k\hat c_n+&\sum_{0\le i\le {k-1}}t_i\la_k\dots\la_{i+1}\hat c_{n_i})
\in X_n\times I^n}$$
where ${\bf t}=(t_1,\ldots,t_k)\in\Delta_k$.
Further details are left to the reader.\qed

\proc{Corollary}Any space has the weak homotopy type of some $\sq$--set.
\qed\endproc

\rk{Notation}We shall often omit the mod signs and use the notation
$C$ both for the $\sq$--set $C$ and its realisation $\mo C\mo$.  We
shall use the full notation whenever there is any possibility of
confusion.

\rk{\label\quad Examples of $\tsq$--sets\key{exxx}}
\let\lrkk\lrk

\lrkk{1}{The singular complex}

Let $X$ be a topological spaces.  The {\sl singular complex} of $X$ is
the $\sq$--set denoted ${\rm S}(X)$ and defined by:

${\rm S}(X)_n=\{f\co I^n\to X\}$,\quad the set of continuous maps and

$\lambda^*(f)=f\circ \lambda $,\quad where $\lambda $ is a face map.

If $g\co X\to Y$ is a continuous map then ${\rm S}(g)\co {\rm S}(X)
\to {\rm S}(Y)$ is the
$\sq$--map defined by ${\rm S}(g)(f)= f\circ g$.  Thus ${\rm S}(.)$ 
is a functor from $Top$ to the category of $\sq$--sets and maps.

An important special case is the singular complex of a point.
This has exactly one cell in each dimension and we shall call
it the {\sl trivial $\sq$--set} and denote it $T$.

\lrkk{2}{The nerve of a category}

(For full details on this example, see [\Trunks].)  Let $I^n_{\rm
cat}$ denote the category whose objects are the vertices of
$I^n$. Morphisms are generated by the oriented edges of $I^n$ and
subject to relations given by taking squares to be commutative
diagrams.  Each face map $\lambda\co I^q\to I^n$ determines a
functor$\lambda_{\rm cat}\co I^q_{\rm cat}\to I^n_{\rm cat}$.  Let
$\C$ be a category; the {\sl nerve} of $\C$, denoted $\N\C$ is the
$\sq$--set defined in analogy to the singular set:

$\N\C_n=\{f\co I^n_{\rm cat}\to \C\}$,\quad the set of functors and

$\lambda^*(f)=f\circ \lambda_{\rm cat} $,\quad where $\lambda $ is a face map.

Given a functor $g\co\C \to\D$ between categories, we have a $\sq$--map
$\N(g)\co\N\C\to\N\D$, given by $\N(g)(f)= f\circ g$.

\lrkk{3}{The rack space}

A {\sl rack} is a set $R$ with a binary operation written $a^b$ such
that $a\mapsto a^b$ is a bijection for all $b\in R$ and such that the
{\sl rack identity}
$$a^{bc}=a^{cb^c}$$ holds for all $a,b,c\in R$.  (Here we use the
conventions for order of operations derived from exponentiation
in arithmetic.  Thus $a^{bc}$ means $(a^b)^c$ and $a^{cb^c}$ means
$a^{c(b^c)}$.)

For example if $a,b\in G$ are two members of a group then conjugation,
given by $a^b:=b\inv ab$, defines a rack structure on $G$.  If $a,b\in \Z$
are two integers then $a^b:=a+1$ defines a rack structure on $\Z$ which
does not come from conjugation in a group.  Another example of a rack
associated to a group  is the {\sl core}
of a group [\Bruck, \Joyce] defined by
$a^b:=ba\inv b$.  For more examples of racks see [\FeRo].

If $R$ is a rack, the {\sl rack space} is the $\sq$--set denoted $BR$ and  
defined by:

$BR_n = R^n$\quad (the $n$--fold cartesian product of $R$ with itself).
$$\d_i^0(x_1,\ldots ,x_n)=(x_1,\ldots
,x_{i-1},x_{i+1},\ldots ,x_n),$$
$$\d_i^1(x_1,\ldots ,x_n)=
((x_1)^{x_i},\ldots ,(x_{i-1})^{x_i},x_{i+1},\cdots
,x_n)
\hbox {\quad for\quad}1\leq i\leq n.$$

More geometrically, we can think of $BR$ as the $\sq$--set with one
vertex, with (oriented) edges labelled by rack elements and with
squares which can be pictured as part of a
link diagram with arcs labelled by $a$, $b$ and $a^b$ (figure 
\figkey\RSpace).

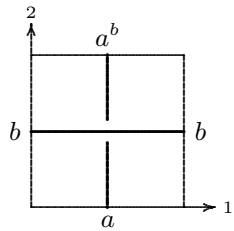
\fig{\RSpace: Diagram of a typical 2--cell of the rack space}
\beginpicture\ninepoint
\setcoordinatesystem units <.4truein,.4truein> point at 0 0
\put {$b$} [r] <-3pt,0pt> at 0 1
\put {$a$} [t] <0pt,-3pt> at 1 0
\put {$a^b$} [b] <0pt,3pt> at 1 2
\put {$b$} [l] <3pt,0pt> at 2 1
\setlinear
\plot 0 0 2 0 2 2 0 2 0 0 /
%\arrow <6pt> [.2, .7] from 1.59 1 to 1.6 1
\arrow <4pt> [.2, .7] from 0 2 to 0 2.4
\arrow <4pt> [.2, .7] from 2 0 to 2.4 0
\setplotsymbol ({\ninerm .})
\plot 1 0 1 0.85 /
\plot 1 1.15 1 2 /
\plot 0 1 2 1 /
\put {$\scriptstyle1$} [l] <2pt,0pt> at 2.4 0
\put {$\scriptstyle2$} [b] <0pt,2pt> at 0 2.4
\endpicture
\endfig

%\eject
%(The orientation in the diagram is determined by the
%axes (labelled 1,2) thus the reflected version represents the
%{\it same} 2--cell.  This point will become important in the next
%section.)

The higher dimensional cubes are determined by the squares in
a way which  was made precise in [\Trunks] and which we now summarise.

A {\sl trunk} $T$ comprises a set $T_0$ of {\sl vertices}, a set $T_1$ of
directed {\sl edges} (between vertices) and a set $T_2$ of oriented squares
of edges (the {\sl preferred squares}).  Thus a trunk is essentially
a cubical set truncated at dimension 2; but notice that a square is
fully determined by its faces (edges).   A trunk $T$ has a {\sl nerve}
$\N T$ which is a $\sq$--set defined in a similar way to the nerve of a
category (above).  The nerve has vertices, edges and squares
in bijection with those of $T$ and higher cubes determined by
the 2--skeleton as follows.  Let $I^n_{\rm trunk}$ be the trunk
determined by the 2--skeleton of $I^n$, then
$\N T_n=\{f\co I^n_{\rm trunk}\to T\}$,  the set of trunk maps,
with face maps given, as usual, by composition.  This construction
is fully compatible with the nerve of a category;  indeed a category
determines a trunk with preferred squares being the commuting
squares;  the two nerves are then identical.

Now a rack $R$ determines the trunk $T(R)$ with one vertex, edges the
set $R$ and preferred squares of the type pictured above.  The rack
spaces $BR$ and the nerve of the trunk $\N T(R)$ coincide.  In
[\Trunks] the trunk formalism is used to establish connections between
the rack space and other classifying spaces.

Notice that the rack space of the rack with one element has
precisely one cube in each dimension.  This is another
description of the trivial $\sq$--set.

\lrkk{4}{The knot space}

Let ${\cal I}_n$ $n=0,1,\ldots$ denote the set of isotopy classes of
smooth immersions (singular knots) $K$ of $S^1$ in $S^3$ with the
following properties.

(1)\qua  The immersion $K$ is an embedding except at $n$ double
points none of which involve the basepoint of $S^1$.

(2)\qua  At each double point the two tangents are linearly independent.

Then ${\cal I}=\{{\cal I}_n\}$ can be made into a $\sq$--set as
follows. Let $K$ be a singular knot defining $[K]\in{\cal
I}_n$. Number the double points of $K$ in order round the knot
starting from the basepoint and define $\d_i^1([K])= [K_i^+]$ where
$K_i^+$ is obtained from $K$ by unwrapping the $i^{\rm th}$ double point to
give a positive crossing and $\d_i^0([K])$ is similarly defined using
a negative crossing.  The face relations given at the start of the
section are easily verified.

\rk{Digression on Vassiliev invariants}

The knot space can be used to interpret Vassiliev invariants.  A {\sl
Vassiliev function} $V=\{V_n\}$ on a $\sq$--set $X$ with values in an
abelian group $G$ is a sequence of functions $V_n\co X_{n}\to G$ for
each $n=0,1,\ldots$ satisfying the {\sl Vassiliev identity}
$$V_n(x)=V_{n-1}(\d^1_ix) -V_{n-1}(\d_i^0x)\quad {\rm for}\quad
i=1,2,\ld,n,\quad x\in X_n.\leqno{\bf\label}$$
If $X$ is the knot space ${\cal I}$ a Vassiliev function is the usual
notion of Vassiliev invariant [\Vass].  The $\sq$--set setting can be
used to interpret the fundamental integration problem for Vassiliev
invariants (see Bar-Natan [\BarNat]) in terms of cohomological
obstructions.  For more details see [\James, section 7].  

\section{James complexes and bundles}

In this section we define, for any $\sq$--set $C$, an infinite family
of associated $\sq$--sets $J^n(C),\>n\ge0$ called the James complexes
of $C$.  These define mock bundles $\zeta^i/C$ which we call James
bundles of $C$.

\rk{Projections}
An {\sl$(n+k,k)$--projection} is a function $\lambda\co I^{n+k}\to I^k$
of the form $$\lambda \co (x_1,x_2,\ldots,x_{n+k})\mapsto (x_{i_1},x_{i_2},
\ldots,x_{i_k}),\leqno{\bf \label}$$ where $1\le i_1<i_2<\ldots<i_k\le n+k$.

Let $P_k^{n+k}$ denote the set of $(n+k,k)$--projections.  Note that
$P_k^{n+k}$ is a set of size $n+k\choose k$.

Let $\lambda\in P_k^{n+k}$ and let $\mu\co I^l\to I^k,\quad l\le k$ be a face
map.  The projection $\mu^\sharp(\lambda)\in P_l^{n+l}$ and the face map
$\mu_\lambda\co I^{n+l}\to I^{n+k}$ are defined uniquely by the 
following pull-back diagram:
$$\matrix{
I^{n+l}&\buildrel\mu_\lambda\over\a&I^{n+k}\cr
\cr
{\scriptstyle\mu^\sharp(\lambda)}\downarrow&&{\scriptstyle\lambda}\downarrow\cr
\cr
I^l&\buildrel\mu\over\a&I^k\cr}
\leqno{\bf \label}$$

\rk{Definition}
Let $C$ be a $\sq$--set.  The {\sl$n^{\rm th}$ associated James
complex} of $C$, denoted $J^n(C)$ is defined as follows.  The
$n$--cells are given by
$$J^n(C)_k=C_{n+k}\times P_k^{n+k}$$
and face maps by
$$\mu^*(x,\lambda)=(\mu_\lambda^*(x),\mu^\sharp(\lambda))$$
where $\mu\co I^l\to I^k,\quad l\le k$ is a face map. 

\rk{Notation}Let $\lambda\in P_k^{n+k}$ and $c\in C_{n+k}$ then
we shall use the notation $c_\lambda$ for the $k$--cube $(c,\lambda)\in
J^n(C)$.  When necessary, we shall use the full notation $(\la_1,\ldots,\la_n)$
for the projection $\lambda$ (given by formula 2.1) where
$\lambda _1<\lambda _2<\ldots<\lambda _n$
 and $\{\lambda _1,\ldots,\lambda _n\}=\{1,\ldots,k+n\}-
\{i_1,\ldots,i_k\}$.  In other words we index cubes of $J^n$ by the
$n$ directions (in order) which are collapsed by the defining
projection.

\sh{Picture for James complexes}

We think of $J^n(C)$ as comprising all the codimension $n$ central
subcubes of cubes of $C$.   For example a 3--cube $c$ of $C$ gives rise
to the three 2--cubes of $J^1(C)$ which are illustrated in
figure \figkey\JaComp.

\fig{\JaComp}
\beginpicture
\setcoordinatesystem units <.8mm,.8mm> point at 0 0
\ninepoint
\put {$\scriptstyle1$} [t] at 7 -0.5
\put {$\scriptstyle2$} [br] at 5 2.5
\put {$\scriptstyle3$} [r] at 0 7
\put {$c_{(2)}$} [tl] at 40 5
\put {$c_{(3)}$} [l] at 50 25
\put {$c_{(1)}$} [b] at 35 41
\putrule from 0 0 to 30 0
\putrule from 0 30 to 30 30
\putrule from 0 0 to 0 30
\putrule from 30 0 to 30 30
\putrule from 20 40 to 50 40
\putrule from 50 10 to 50 40
\plot 0 30  20 40 /
\plot 30 30  50 40 /
\plot 30 0  50 10 /
\setdashes
\putrule from 20 10  to 20 40
\putrule from 20 10  to 50 10
\plot 0 0 20 10 /
\setsolid
\linethickness=1.2pt
\putrule from 0 15 to 30 15
\putrule from 15 0 to 15 30
\putrule from 10 35 to 40 35
\putrule from 40 5 to 40 35
\putrule from 10 5 to 10 35
\putrule from 10 5 to 40 5
\putrule from 35 10 to 35 40 
\putrule from 20 25 to 50 25
\putrule from 25 5 to 25 35
\putrule from 10 20 to 40 20
\setlinear\setplotsymbol ({\tenrm .})
\plot 15 30  35 40 /
\plot 30 15  50 25 /
\plot 15 0  35 10 /
\plot 15 15  35 25 /
\plot 0 15  20 25 /
\endpicture
\endfig

In figure \JaComp\ we have used full notation for projections. 
Thus for example,
$c_{(2)}$ corresponds to the projection $(x_1,x_2,x_3)\mapsto(x_1,x_3)$,
($x_2$ being collasped).

This picture can be made more precise by considering the section
$s_\lambda\co I^k\to I^{n+k}$ of $\lambda$ given by
$$
\def\ha{{\scriptstyle{1\over2}}}
s_\lambda(x_1,x_2,\ldots,x_k)=(\ha,\ldots,\ha,x_1,\ha,\ldots,\ha,x_2,\ha,
\ldots,\ha,x_k,\ha,\ldots)
$$
where the non-constant coordinates are in places $i_1,i_2,\ldots,i_k$
and $\lambda$ is given by 2.1.

For example in the picture the image of $s_\lambda$ where 
$\lambda(x_1,x_2,x_3)\mapsto(x_1,x_3)$ is the 2--cube labelled $c_{(2)}$.

Now the commuting diagram (2.2) which defines the face maps implies
that the $s_\lambda$'s are compatible with faces and hence they fit 
together to define a map 
$$p_n\co\mo J^n(C)\mo\to\mo C\mo\quad\hbox{given by}\quad 
p_n[c_\lambda,t]=[c,s_\lambda(t)].$$  

\sh{James bundles}

We now observe that $p_n\co\mo J^n(C)\mo\to\mo C\mo$
is a mock bundle projection.  We shall need to define a mock bundle
with base a $\sq$--set.  Mock bundles with base a cell complex are
defined (in a PL setting) in [\BRS] and in [\Rackspace] we extend the 
concept to define a smooth mock bundle with base a smooth CW complex.

For the purposes of this paper, we shall not need to pay attention to 
smooth technicalities and shall only need a minor extension of the 
treatment in [\BRS].

\proclaim{Definition}Mock bundle over a $\sq$--set.\rm

Let $C$ be a $\sq$--set.   A {\sl mock bundle $\xi$ over 
$C$} of codimension $q$ (denoted $\xi^q/C$) comprises
a {\sl total space} $E_\xi$ and a {\sl projection} 
$p_\xi\co E_\xi\to \mo C\mo$ with the following property.%
\fnote{Note that the notation used here for dimension of
a mock bundle, namely that $q$ is {\it codimension}, is the negative of that used
in [\BRS] where $\xi^q/C$ meant a mock bundle of {\it fibre dimension}
$q$ ie codimension $-q$.  The notation used here is consistent with the
usual convention for cohomology.}

Let $c$ be an $n$--cell of $\mo C\mo$ with characteristic map 
$\chi_c\co I^n\to \mo C\mo$, then there is a manifold (with corners)
$B_c$ of dimension $n-q$ called the {\sl block} over $c$
and a proper map $p_c\co B_c\to I^n$ and a map $b_c\co B_c\to E_\xi$
such that the following diagram is a pull-back:
$$
\matrix{
B_c&\buildrel b_c\over \a&E_\xi\cr
&&\cr
\downarrow\!{\scriptstyle p_c}&&\downarrow\!{\scriptstyle p_\xi}\cr
&&\cr
I^n&\buildrel \chi_c\over \a&\mo C\mo\cr}
$$

\rk{Key example}{\sl James bundles}

Consider the projection
$$p_n\co \mo J^n(C)\mo\to\mo C\mo.$$  
Where $J^n(C)$ is the $n^{\rm th}$ associated James complex of the $\sq$--set $C$.

If we choose
a particular $(n+k)$--cell $\sigma$ of $C$ then the pull back of $p_n$
over $I^{n+k}$ (by the characteristic map for $\sigma$) is
a $k$--manifold (in fact it is the $n+k\choose k$ copies of 
$I^k$ corresponding to the elements of $P^{n+k}_k$).  Therefore
$p_n$ is the projection of a mock bundle of codimension
$n$, which we shall call the {\sl$n^{\rm th}$ James bundle} of $C$
denoted $\zeta^n(C)$.  

\rk{Embedding the James bundles in $\mo C\mo\times\re$}

Let $C$ be a $\sq$--set.  The James bundles can be embedded
in $\mo C\mo\times\re$.  This is done by ordering
the cubes of $J^n(C)$ over a particular cube of $C$ and 
lifting in that order.  Recall that the $k$--cubes of $J^n(C)$
lying over a $(k+n)$--cube are indexed by projections 
$\lambda=(i_1,\ldots,i_n)\in P^{n+k}_n$.  These may
be ordered lexicographically. The lexicographic order is compatible 
with face maps and can be used to define the required embedding
by induction on dimension of cells of $C$ as follows.  

Suppose inductively that the embedding has been defined 
over cells of $C$ of dimension $\le k+n-1$.  

Consider a $(k+n)$--cube $c\in C$ with characteristic map
$\chi_c\co I^{k+n}\to \mo C\mo$.  Pulling the embedding
back (where it is already defined) over $\chi_c$ gives an
embedding of $\zeta^n(\d I^{k+n})$ in $\d I^{k+n}\times\re$.
Now embed the centres of the
$k$--cubes of $J^n(I^{k+n})$ at $(\ha,\ldots,\ha)\times r_\lambda$
where $r_\lambda$ are real numbers for $\lambda\in P^{n+k}_n$ chosen
to increase strictly corresponding to the lexicographic order
on $P^{n+k}_n$.  

Now embed each $k$--cube of $J^n(I^{k+n})$ as the cone
on its (already embedded) boundary.  The resulting may be
smoothed if desired and then pushed forwards to
$\mo C\mo\times\re$ using $\chi_c\times {\rm id}$.

In the next section we shall give precise smooth formul\ae\ for this
embedding using a bump function.

The embedding is in fact framed.  This can be seen as follows.
Each $k$--cube $p_n(c,\la)$, where $\lambda=(i_1,\ldots,i_n)$, of
$J^n(I^{k+n})$ is framed in $I^{k+n}$ by the $n$ vectors parallel
to directions $i_1,\ldots,i_n$.  These lift to parallel vectors
in $I^{n+k}\times\re$ and the framing is completed by the vector
parallel to the positive $\re$ direction (vertically up).  This
framing is compatible with faces and defines a framing of $\zeta^n(C)$
in $\mo C\mo\times\re$.  The formul\ae\ which we shall give in the next
section also give formul\ae\ for the framing.

For the special case $n=1$ the map of $\zeta^n(C)$ to $\re$ can be
simply described: the centre of $c_{(k)}$ is mapped to $k$. This
determines a map, linear on simplexes of $Sd_\Delta \zeta^1(C)$, to
$\re$. It follows that the centre of $\d_i^\epsilon c_{(k)}$ is mapped
to $k$ if $i\ge k$ and to $k-1$ if $i<k$.

In figure \JaComp\ we illustrated $J^1(C)$  for a 3--cube
$c\in C$.  The embeddings in $\mo C\mo\times\re$ (before smoothing) above
each of the three 2--cubes are illustrated in the figure \figkey\JaEmb.
\fig{\JaEmb}
\beginpicture
\setcoordinatesystem units < 0.7cm, 0.7cm>
\unitlength= 0.7cm
\linethickness=1pt
\setplotsymbol ({\makebox(0,0)[l]{\tencirc\symbol{'160}}})
\setshadesymbol ({\thinlinefont .})
\setlinear
%
% Fig POLYLINE object
%
\linethickness= 0.500pt
\setplotsymbol ({\thinlinefont .})
\plot  0.921 20.034  2.826 20.415 /
\plot  2.826 20.415  4.731 20.034 /
\plot  4.731 20.034  4.921 21.654 /
\plot  4.921 21.654  5.683 22.892 /
\plot  5.683 22.892  3.778 23.273 /
\plot  3.778 23.273  1.873 22.892 /
\plot  1.873 22.892  1.111 21.654 /
\plot  1.111 21.654  0.921 20.034 /
\putrule from  0.921 20.034 to  0.921 20.034
%
% Fig POLYLINE object
%
\linethickness= 0.500pt
\setplotsymbol ({\thinlinefont .})
\plot  3.112 22.415  1.111 21.654 /
%
% Fig POLYLINE object
%
\linethickness= 0.500pt
\setplotsymbol ({\thinlinefont .})
\plot  3.112 22.415  1.873 22.892 /
%
% Fig POLYLINE object
%
\linethickness= 0.500pt
\setplotsymbol ({\thinlinefont .})
\plot  3.112 22.415  3.778 23.273 /
%
% Fig POLYLINE object
%
\linethickness= 0.500pt
\setplotsymbol ({\thinlinefont .})
\plot  3.112 22.415  5.683 22.892 /
%
% Fig POLYLINE object
%
\linethickness= 0.500pt
\setplotsymbol ({\thinlinefont .})
\plot  3.112 22.415  4.921 21.654 /
%
% Fig POLYLINE object
%
\linethickness= 0.500pt
\setplotsymbol ({\thinlinefont .})
\plot  3.112 22.415  2.826 20.415 /
\plot  3.112 22.415  4.731 20.034 /   %%%%
\plot  3.112 22.415  0.921 20.034 /   %%%%
%
% Fig POLYLINE object
%
\linethickness= 0.500pt
\setplotsymbol ({\thinlinefont .})
\putrule from  0.921 20.034 to  4.731 20.034
\plot  4.731 20.034  5.683 22.892 /
\putrule from  5.683 22.892 to  3.683 22.892
%
% Fig POLYLINE object
%
\linethickness= 0.500pt
\setplotsymbol ({\thinlinefont .})
\plot  0.921 20.034  1.492 21.654 /
%
% Fig POLYLINE object
%
\linethickness= 0.500pt
\setplotsymbol ({\thinlinefont .})
\plot  1.873 22.892  1.587 21.939 /
%
% Fig POLYLINE object
%
\linethickness= 0.500pt
\setplotsymbol ({\thinlinefont .})
\putrule from  1.873 22.892 to  3.302 22.892
%
% Fig POLYLINE object
%
\linethickness= 0.500pt
\setplotsymbol ({\thinlinefont .})
\putrule from 12.351 20.034 to 16.161 20.034
\plot 16.161 20.034 17.113 22.892 /
\putrule from 17.113 22.892 to 13.208 22.892
\plot 13.208 22.892 12.351 20.034 /
%
% Fig POLYLINE object
%
\linethickness= 0.500pt
\setplotsymbol ({\thinlinefont .})
\plot  6.636 20.034  8.636 20.415 /
\plot  8.636 20.415 10.446 20.034 /
%
% Fig POLYLINE object
%
\linethickness= 0.500pt
\setplotsymbol ({\thinlinefont .})
\plot  7.588 22.892  9.493 23.273 /
%
% Fig POLYLINE object
%
\linethickness= 0.500pt
\setplotsymbol ({\thinlinefont .})
\plot  9.493 23.273 11.398 22.892 /
%
% Fig POLYLINE object
%
\linethickness= 0.500pt
\setplotsymbol ({\thinlinefont .})
\plot  9.493 23.273  8.636 20.415 /
%
% Fig POLYLINE object
%
\linethickness= 0.500pt
\setplotsymbol ({\thinlinefont .})
\putrule from  9.684 22.892 to 11.398 22.892
\plot 11.398 22.892 10.446 20.034 /
\putrule from 10.446 20.034 to  6.636 20.034
\plot  6.636 20.034  7.588 22.892 /
\putrule from  7.588 22.892 to  9.112 22.892
%
% Fig TEXT object
%
\def\SetFigFont{\ninepoint}%
\put{\SetFigFont $c(2)$ } [lB] <0pt, 10pt> at  8.446 18.605
%
% Fig TEXT object
%
\put{\SetFigFont $c(3)$} [lB] <0pt, 10pt> at  2.445 18.605
%
% Fig TEXT object
%
\put{\SetFigFont $c(1)$} [lB] <0pt, 10pt> at 14.161 18.701
%
% Fig TEXT object
%
\put{\SetFigFont$\scriptstyle 1$} [lB] <0pt,2pt> at  1.302 19.558
%
% Fig TEXT object
%
\put{\SetFigFont$\scriptstyle  2$} [lB] <0pt, 5pt> at  1.302 20.701
%
% Fig TEXT object
%
\put{\SetFigFont$\scriptstyle  2$} [lB] <0pt, 0pt> at 12.732 19.653
%
% Fig TEXT object
%
\put{\SetFigFont$\scriptstyle  3$} [lB] <-4pt, 0pt> at 12.827 20.701
%
% Fig TEXT object
%
\put{\SetFigFont$\scriptstyle  1$} [lB] at  7.017 19.653
%
% Fig TEXT object
%
\put{\SetFigFont$\scriptstyle  3$} [lB] <-2pt, 0pt> at  7.017 20.606
\linethickness=0pt
\putrectangle corners at  0.921 23.273 and 17.113 18.605
\endpicture
\endfig

\section{James--Hopf invariants}

In this section we shall show that James complexes are strongly
connected with classical James--Hopf invariants; indeed they define
(generalised) James--Hopf invariants for any $\sq$--set $C$ and this
explains our choice of terminology.  We do this by looking carefully
at James' original contruction [\JRed, \JHopf] and in the process we
give precise formul\ae\ for the framed embedding of the James bundles
in $\mo C\mo\times\re$.  

\sh{The James construction}

\rk{Definition}{\sl Free topological monoid}

Let $X$ be a topological space based at $*$.  The {\sl free
topological monoid} on $X$, denoted $X_\infty$ comprises all words
$x_1\cdots x_m$, where $x_i\in X$ and $m\ge 0$ with the
identifications given by regarding $*$ as a unit; ie $x_1\cdots
x_{i-1}*x_{i+1}\cdots x_m\sim x_1\cdots x_{i-1}x_{i+1}\cdots x_m$ for
$i=1,\ldots,m$. Thus $X_\infty$ is a quotient of the disjoint union
$\coprod_{m\ge0}X^m$ and this defines the topology. The space
$X_\infty$ is based at the empty word which we also denote by $*$.  
A point of $X_\infty$ has a unique representation as a {\sl reduced}
word, ie one with $x_i\ne*$ for all $i$.

\rk{The James map}

Let $X$ be a based topological space and suppose that there is
a map $\rho_X\co X\to I$ with the property that $\rho_X^{-1}(0)=\{*\}$.
Let $\Omega(X)$ denote the loop space of $X$ and $S(X)=X\wedge S^1$ 
the suspension of $X$.  (Notice that we take the suspension
coordinate second.)
We shall identify $S^1$ with $I/\d I$. Thus a point 
of $S^1-\{*\}$ is uniquely represented as a real number $t$ with $0<t<1$,
and a point of $S(X)-\{*\}= (X\wedge S^1)-\{*\}$ is uniquely represented as
a pair $(x,t)$ where $x\in X-\{*\}$ and $0<t<1$.  It is convenient
to regard the suspension coordinate as vertical.

Suppose given a non-empty reduced word $x_1\cdots x_m\in X_\infty$
where each $x_i\ne *$.  Let $\alpha_i=\rho_X(x_i)$, and let $t_i$ be
given by $t_0=0$ and
$$t_i={\alpha_1+\dots +\alpha_i\over\alpha_1+\dots +\alpha_m}$$
for $0< i\le m$.  Then $0=t_0<t_1<\ldots<t_m=1$.

The {\sl James map} $k_X\co X_\infty\to\Omega S(X)$
is given by $k_X(*)=*$ and for non-empty words
$$k_X(x_1\cdots x_m)(t)=\Bigl(x_i,{t-t_{i-1}\over t_i-t_{i-1}}
\Bigr),\leqno{\bf \label}$$
for $t_{i-1}\le t\le t_i$.

In words, what $k_X$ does is to map the word $x_1\cdots x_m$
to a loop in $S(X)$ which comprises $m$ vertical loops passing
through $x_1, \dots, x_m$ respectively, with the time parameters
adjusted by using $\rho$ to make the time spent on a subloop
go to zero as the corresponding point $x_i$ moves to the basepoint
of $X$.

James proves that $k_X$ is a homotopy equivalence if $X$ is
a countable CW complex with just one vertex.  However the 
result extends to a considerably more
general class of spaces, [\Puppe].  We shall be particularly interested
in the case when $X$ is the $n$--sphere and we shall abbreviate
$k_{S^n}$ to $k_n$.

\rk{The classical James--Hopf invariants}

We make the following further identifications:
$$S^{n+1}=S(S^n)=S^n\wedge S^1=S^1\wedge\cdots\wedge S^1=
I/\d I\wedge\cdots\wedge I/\d I=I^{n+1}/\d I^{n+1}.$$
This means that a point of $S^n-\{*\}$ is uniquely
an $n$-tuple $x=(x_1,x_2,\ldots,x_n)$ where $0<x_i<1$. 

We choose the map $\rho_{S^n}$ rather carefully.

Let $\rho\co I\to I$ be a smooth bump function with
the property that $\rho^{-1}(0)=\{0,1\}$
and all derivatives vanish at $0$ and $1$. Let $\rho_{S^1}\co S^1\to I$
be the induced map, and define $\rho_{S^n}\co S^n\to I$
by $\rho_{S^n}(x_1,x_2,\ldots,x_n)=\rho(x_1)\rho(x_2)\dots\rho(x_n)$.

Now let $g_n\co S^1_\infty\to S^n_\infty$ be given by
$$g_n(x_1\cdots x_m)=
\prod_\la (x_{\la_1},\ldots,x_{\la _n}),\leqno{\bf\label}$$
where the product is over all strictly monotone $\la\co\{1,\ldots,n\}\to\{1,\ldots,m\}$
and is taken in lexicographical order.

The $n$--{\sl th James--Hopf invariant} $J_n\co [S(X),S^2]\to
[S(X),S^{n+1}]$ is
induced by the composition
$$\Omega S^2{\buildrel \tau\over\a}S^1_\infty{\buildrel g_n\over\a}
S^n_\infty{\buildrel k_n\over\a}\Omega S^{n+1},$$
where $\tau$ is a homotopy inverse to $k_1$.

\sh{The trivial $\tsq$--set $T$}

We now turn to the connection between embedded James
bundles and the James--Hopf invariants. We start by
identifying the homotopy type of the trivial
$\sq$--set $T$
(with exactly one cell in each dimension).

\proc{Proposition}The realisation of the $\sq$--set $T$ can be identified
with $S^1_\infty$ and hence has the homotopy type of $\Omega(S^2)$.

\prf Denote the unique $m$--cube of $T$ by $c_m$.   Let $x\in\mo T\mo$
then $x\in c_m^\circ$ for some $m$ and hence $x$ has a unique expression 
as $(x_1,\ldots,x_m)$ where each $0<x_i<1$.  Moreover the glue given
by the face maps has the effect of omitting $x_i$'s which become 0
or 1.  It follows that the map $x=(x_1,\ldots,x_m)\mapsto x_1\cdots x_m$
is a homeomorphism $\mo T\mo\to S^1_\infty$. \qed

\rk{Remark}Antolini [\Ant] has generalised this result and 
shown that the singular $\sq$--set
of any topological space $X$ has the weak homotopy type of 
$X\times \Omega(S^2)$, see also [\AntWie]. 

\sh{The James bundles of $T$}

Now consider the $n^{\rm th}$ James bundle $\zeta^n(T)$.  In the last section
we showed that it embeds as a framed mock bundle in $\mo T\mo\times\re$.
If we apply the Thom--Pontrjagin construction to this framed embedding
we obtain a map $\mo T\mo\times\re\to S^{n+1}$ and hence a map
$q_n\co S(\mo T\mo)\to S^{n+1}$.  Using the last proposition we
have a map
$$q_n\co S(S^1_\infty)\a S^{n+1}.$$

The following proposition connects the James bundles of $T$ with
the James--Hopf invariants:

\proc{Proposition}The framed embedding of $E_{\zeta^n(T)}$ in 
$\mo T\mo\times\re$
can be chosen so that $q_n$ is the adjoint of
$$S^1_\infty{\buildrel g_n\over\a}S^n_\infty{\buildrel k_n\over\a}
\Omega(S^{n+1}).$$

\prf Let $h_n\co I^{n+k}\times I\to S^{n+1}$ be the map determined
by the adjoint of  $k_ng_n$ restricted to the $(n+k)$--cell of $T$. Then
$h_n$ is transverse to $a=(\ha,\ldots,\ha)$ and we shall see
that the framed submanifold
$h^{-1}(a)$ can be identified with the block of 
$E_{\zeta^n(T)}$ over $I^{n+k}$.  Moreover
the restriction of $h_n$ to $h_n^{-1}({\rm int}(I^{n+1}))$ 
provides a trivialisation
of an open tubular neighbourhood of  
$E_{\zeta^n(T)}$ in $I^{n+k}\times I\subset I^{n+k}\times\re$.
The closure of this neighbourhood is the whole of $I^{n+k}\times I$. 
Figure \figkey\Frames\ illustrates the case $n=k=1$. Only the vertical part 
of the framing is shown for clarity. In general the $k$--cubes 
$M_{c_\la}$ are placed in lexicographical
order above their images $s_\la(I^k)$ in $I^{n+k}$.

\fig{\Frames}
\beginpicture
\setcoordinatesystem units < 1.000cm, 1.000cm>
\unitlength= 1.000cm
\linethickness=1pt
\setplotsymbol ({\makebox(0,0)[l]{\tencirc\symbol{'160}}})
\setshadesymbol ({\thinlinefont .})
\setlinear
%
% Fig POLYLINE object
%
\linethickness= 0.500pt
\setplotsymbol ({\thinlinefont .})
\plot  3.937 24.320  3.016 22.892 /
\putrule from  3.016 22.892 to  0.889 22.892
%
% Fig POLYLINE object
%
\linethickness= 0.500pt
\setplotsymbol ({\thinlinefont .})
\plot  0.889 24.797  1.841 26.194 /
\plot  1.841 26.194  3.937 26.226 /
%
% Fig POLYLINE object
%
\linethickness= 0.500pt
\setplotsymbol ({\thinlinefont .})
\putrule from  0.889 24.797 to  0.889 22.892
%
% Fig POLYLINE object
%
\linethickness= 0.500pt
\setplotsymbol ({\thinlinefont .})
\putrule from  3.016 24.797 to  3.016 22.892
%
% Fig POLYLINE object
%
\linethickness= 0.500pt
\setplotsymbol ({\thinlinefont .})
\putrule from  3.937 26.226 to  3.937 24.320
%
% Fig POLYLINE object
%
\linethickness= 0.500pt
\setplotsymbol ({\thinlinefont .})
\plot  2.000 22.892  2.921 24.320 /
%
% Fig POLYLINE object
%
\linethickness= 0.500pt
\setplotsymbol ({\thinlinefont .})
\plot  3.937 26.226  3.016 24.797 /
\putrule from  3.016 24.797 to  0.889 24.797
%
% Fig POLYLINE object
%
\linethickness= 0.500pt
\setplotsymbol ({\thinlinefont .})
\setdots < 0.1270cm>
\plot  1.397 23.654  0.889 22.892 /
%
% Fig POLYLINE object
%
\linethickness= 0.500pt
\setplotsymbol ({\thinlinefont .})
\setsolid
\putrule from  1.492 25.527 to  1.492 23.654
%
% Fig POLYLINE object
%
\linethickness= 0.500pt
\setplotsymbol ({\thinlinefont .})
\putrule from  1.715 25.527 to  1.715 23.749
%
% Fig POLYLINE object
%
\linethickness= 0.500pt
\setplotsymbol ({\thinlinefont .})
\putrule from  3.397 25.527 to  3.397 23.654
%
% Fig POLYLINE object
%
\linethickness= 0.500pt
\setplotsymbol ({\thinlinefont .})
\putrule from  2.318 24.448 to  2.318 23.368
%
% Fig POLYLINE object
%
\linethickness= 0.500pt
\setplotsymbol ({\thinlinefont .})
\putrule from  2.413 24.448 to  2.413 23.527
%
% Fig POLYLINE object
%
\linethickness= 0.500pt
\setplotsymbol ({\thinlinefont .})
\plot  2.508 24.416  2.508 24.416 /
%
% Fig POLYLINE object
%
\linethickness= 0.500pt
\setplotsymbol ({\thinlinefont .})
\putrule from  2.508 24.384 to  2.508 23.685
%
% Fig POLYLINE object
%
\linethickness= 0.500pt
\setplotsymbol ({\thinlinefont .})
\putrule from  2.603 24.320 to  2.603 23.844
%
% Fig POLYLINE object
%
\linethickness= 0.500pt
\setplotsymbol ({\thinlinefont .})
\putrule from  1.905 25.527 to  1.905 24.130
%
% Fig POLYLINE object
%
\linethickness= 0.500pt
\setplotsymbol ({\thinlinefont .})
\putrule from  2.000 25.527 to  2.000 24.892
%
% Fig POLYLINE object
%
\linethickness= 0.500pt
\setplotsymbol ({\thinlinefont .})
\putrule from  2.127 25.527 to  2.127 24.924
%
% Fig POLYLINE object
%
\linethickness= 0.500pt
\setplotsymbol ({\thinlinefont .})
\putrule from  2.254 25.527 to  2.254 24.575
%
% Fig POLYLINE object
%
\linethickness= 0.500pt
\setplotsymbol ({\thinlinefont .})
\putrule from  2.381 25.527 to  2.381 24.543
%
% Fig POLYLINE object
%
\linethickness= 0.500pt
\setplotsymbol ({\thinlinefont .})
\putrule from  2.508 25.527 to  2.508 24.511
%
% Fig POLYLINE object
%
\linethickness= 0.500pt
\setplotsymbol ({\thinlinefont .})
\putrule from  3.080 25.527 to  3.080 23.876
%
% Fig POLYLINE object
%
\linethickness= 0.500pt
\setplotsymbol ({\thinlinefont .})
\putrule from  2.953 25.495 to  2.953 24.066
%
% Fig POLYLINE object
%
\linethickness= 0.500pt
\setplotsymbol ({\thinlinefont .})
\putrule from  2.603 25.527 to  2.603 24.448
%
% Fig POLYLINE object
%
\linethickness= 0.500pt
\setplotsymbol ({\thinlinefont .})
\putrule from  1.587 25.527 to  1.587 23.654
%
% Fig POLYLINE object
%
\linethickness= 0.500pt
\setplotsymbol ({\thinlinefont .})
\putrule from  3.493 23.622 to  3.493 25.527
\putrule from  3.493 25.527 to  1.397 25.527
%
% Fig POLYLINE object
%
\linethickness= 0.500pt
\setplotsymbol ({\thinlinefont .})
\putrule from  3.175 25.527 to  3.175 23.717
%
% Fig POLYLINE object
%
\linethickness= 0.500pt
\setplotsymbol ({\thinlinefont .})
\putrule from  3.270 25.527 to  3.270 23.654
%
% Fig POLYLINE object
%
\linethickness= 0.500pt
\setplotsymbol ({\thinlinefont .})
\putrule from  2.699 24.257 to  2.699 24.003
%
% Fig POLYLINE object
%
\linethickness= 0.500pt
\setplotsymbol ({\thinlinefont .})
\putrule from  2.699 25.749 to  2.699 24.416
%
% Fig POLYLINE object
%
\linethickness= 0.500pt
\setplotsymbol ({\thinlinefont .})
\putrule from  1.397 25.463 to  1.397 23.622
%
% Fig POLYLINE object
%
\linethickness= 0.500pt
\setplotsymbol ({\thinlinefont .})
\putrule from  2.889 26.194 to  2.889 24.289
%
% Fig POLYLINE object
%
\linethickness= 0.500pt
\setplotsymbol ({\thinlinefont .})
\putrule from  2.191 24.543 to  2.191 23.209
%
% Fig POLYLINE object
%
\linethickness= 0.500pt
\setplotsymbol ({\thinlinefont .})
\putrule from  2.794 26.067 to  2.794 24.257
%
% Fig POLYLINE object
%
\linethickness= 0.500pt
\setplotsymbol ({\thinlinefont .})
\putrule from  2.095 23.050 to  2.095 24.829
%
% Fig POLYLINE object
%
\linethickness= 0.500pt
\setplotsymbol ({\thinlinefont .})
\putrule from  1.810 25.495 to  1.810 23.971
%
% Fig POLYLINE object
%
\linethickness= 0.500pt
\setplotsymbol ({\thinlinefont .})
\putrule from  2.000 24.797 to  2.000 22.892
\linethickness= 0.500pt
\setplotsymbol ({\thinlinefont .})
\setdots < 0.0635cm>
%
% Fig INTERPOLATED PT SPLINE
%
\plot  1.969 24.797 	 2.095 24.860
	 2.148 24.771
	 2.159 24.670
	 2.168 24.606
	 2.191 24.543
	 2.277 24.486
	 2.381 24.479
	 2.479 24.545
	 2.551 24.647
	 2.602 24.761
	 2.635 24.860
	 2.647 24.953
	 2.645 25.069
	 2.638 25.183
	 2.635 25.273
	 2.639 25.371
	 2.646 25.495
	 2.656 25.620
	 2.667 25.718
	 2.678 25.774
	 2.695 25.845
	 2.730 25.971
	 2.761 26.050
	 2.788 26.111
	 2.826 26.194
	/
\linethickness= 0.500pt
\setplotsymbol ({\thinlinefont .})
%
% Fig INTERPOLATED PT SPLINE
%
\plot  1.397 23.654 	 1.491 23.652
	 1.561 23.657
	 1.651 23.685
	 1.729 23.769
	 1.793 23.887
	 1.850 24.008
	 1.905 24.098
	 1.990 24.203
	 2.101 24.327
	 2.165 24.385
	 2.233 24.433
	 2.305 24.466
	 2.381 24.479
	 2.470 24.468
	 2.555 24.434
	 2.637 24.383
	 2.713 24.320
	 2.784 24.251
	 2.848 24.182
	 2.953 24.066
	 3.009 23.983
	 3.067 23.872
	 3.132 23.762
	 3.207 23.685
	 3.308 23.655
	 3.386 23.651
	 3.493 23.654
	/
\linethickness=1pt
\setplotsymbol ({\makebox(0,0)[l]{\tencirc\symbol{'160}}})
\setdashes < 0.1429cm>
%
% Fig INTERPOLATED PT SPLINE
%
\plot  1.397 24.479 	 1.489 24.590
	 1.569 24.684
	 1.641 24.763
	 1.704 24.829
	 1.812 24.926
	 1.905 24.987
	 2.016 25.028
	 2.084 25.045
	 2.156 25.058
	 2.228 25.069
	 2.298 25.077
	 2.413 25.082
	 2.472 25.079
	 2.543 25.070
	 2.622 25.058
	 2.703 25.041
	 2.784 25.022
	 2.861 25.001
	 2.985 24.955
	 3.071 24.905
	 3.176 24.826
	 3.240 24.773
	 3.313 24.709
	 3.396 24.632
	 3.493 24.543
	/
\linethickness=1pt
\setplotsymbol ({\makebox(0,0)[l]{\tencirc\symbol{'160}}})
%
% Fig INTERPOLATED PT SPLINE
%
\plot  2.000 23.876 	 2.122 23.929
	 2.211 23.970
	 2.318 24.035
	 2.391 24.105
	 2.472 24.200
	 2.548 24.300
	 2.603 24.384
	 2.646 24.473
	 2.691 24.588
	 2.732 24.704
	 2.762 24.797
	 2.779 24.857
	 2.799 24.939
	 2.824 25.054
	 2.840 25.126
	 2.857 25.210
	/
\linethickness=0pt
\putrectangle corners at  0.889 26.226 and  3.937 22.892
\endpicture
\endfig

We now give a precise formula for the embedding of this
block and then the statements made above can all be checked.
Let $\la^1,\dots,\la^N$ be the list of $\la$'s in the
definition of $g_n$ (3.2) with $m=n+k$,  taken in order. Then
from 3.1 and 3.2 we have
$$h_n(x_1,\ldots,x_{n+k},t)=
\Bigl(x_{\la_1^i},\ld,x_{\la^i_n},{t-t_{i-1}\over t_i-t_{i-1}}
\Bigr),$$
$$\hbox{where}\quad
t_i={\alpha_1+\dots +\alpha_i\over\alpha_1+\dots +\alpha_N}\quad
\hbox{for}\quad 0\le i\le N\quad\hbox{and}\quad
\alpha_i=\rho_n(x_{\la_1^i},\ld,x_{\la^i_n}).$$
Each $\la^i=(\la^i_1,\ldots,\la^i_n)$ corresponds to a $k$--cube
$M_{c_{\la^i}}$ of the $n^{\rm th}$ associated James complex, namely the cube
$(c_{n+k},\la^i)\in J^n(T)_k$, in the block of $E_{\zeta^n(T)}$ over 
$c_{n+k}$. Recall that the $\la^i_j$'s index the directions
in $I^{n+k}$ normal to the $k$--cube $(c_{n+k},\la^i)$.

We will define an embedding
$$e_{\la^i}\co I^k\times {\rm int}(I^n)\times {\rm int}(I)\to 
I^{n+k}\times I\subset I^{n+k}\times\re.$$

Let $sh_{\la^i}\co I^{n+k}\to I^{n+k}$ be the isometry given by permuting
coordinates as follows: $sh_{\la^i}$ sends the
 $(k+s)^{\rm th}$ coordinate to the ${\la^i_s}^{\rm th}$ coordinate,
 for $1\le s\le n$, and preserves order
on the remaining $k$ coordinates.  Write $sh_{\la^i}(y)=x$, then
$$e_{\la^i}(y,u)=(x,ut_i+(1-u)t_{i-1}).$$
Note that the images of the $e_{\la^i}$ are disjoint.

Now define the embedding of the $k$--cube $(c_{n+k},\la^i)\in J^n(T)_k$
by identifying $I^k$ with $I^k\times(\ha,\ld,\ha)\times\ha
\subset I^k\times {\rm int}(I^n)\times {\rm int}(I)$ and composing
$e_{\la^i}$ with $\chi_{c_{n+k}}\times{\rm id}_\re$.  If we do
this for each of the $\la^i$ we obtain
a framed embedding of the block of $\zeta^n(T)$
over $c_{n+k}$ and the required properties follow from comparing
the formul\ae\ above with 3.1 and 3.2.
\qed

\rk{Remark}
The embeddings $e_{\la^i}$ above extend to maps (not embeddings)
of $I^k\times I^n\times I\to I^{n+k}\times \re$ whose images fill the
whole of $I^{n+k}\times I$.  However the Thom--Pontrjagin
construction works for such framings, since the boundary of
the tubular neighbourhood is mapped to the basepoint.  If a closed
trivial tubular neighbourhood is required, then we can
restrict to a small $(k+1)$--cube at the centre of $I^k\times I$.

We now turn to general $\sq$--sets.  The following proposition
follows at once from definitions.

\proc{Proposition}James bundles are natural, ie given a $\sq$--map
$f\co C\to D$ then $\zeta^n(C)$ is the pull-back $f^*(\zeta^n(D))$.
\qed

\sh{The James--Hopf invariants of a $\tsq$--set}

If $C$ is any $\sq$--set there is a canonical {\sl constant}
$\sq$--map $t_C\co C\to T$, so proposition 3.5 implies that
$\zeta^n(C)=t_C^*\zeta^n(T)$. 
Further the framed embedding of $\zeta^n(T)$ in $\mo T\mo\times\re$
pulls back to a framed embedding of $\zeta^n(C)$ in $\mo C\mo\times\re$
which (by the Thom--Pontrjagin construction) determines a map
$S(C)\to S^{n+1}$.  

We call the homotopy class of this map 
$j_n(C)\in[S(C),S^{n+1}]$
the {\sl $n^{\rm th}$ James--Hopf invariant} of the $\sq$--set $C$.
%\eject
\proc{Theorem}

{\rm(1)}\qua The James--Hopf invariants of a $\sq$--set are natural, ie
given a $\sq$--map $f\co C\to D$ then $j_n(C)=S(f)^*j_n(D)$.

{\rm(2)}\qua If we
identify $T$ with $S^1_\infty$ by proposition 3.3 and with
$\Omega(S^2)$ via the homotopy equivalence $k_1$ then the 
classical James--Hopf invariant  $J_n$ (defined earlier)
is induced by composition with $j_n(T)$.

{\rm(3)}\qua Let $q_C\co S(C)\to S^2$ denote the adjoint of 
$C{\buildrel t_C\over\a}T=S^1_\infty{\buildrel k_1\over\a}
\Omega(S^2)$ then $j_n(C)=J_n(q_C)$.

\prf  Naturality follows at once from definitions.  Moreover
(2) is a restatement of proposition
3.4 after observing that $q_n\equiv j_1(T)$
(the notation $q_n$ was used in 3.4).  Finally
$$J_n(q_C)={\rm adj}((k_n g_n)t_C)=S(t_C)^*{\rm adj}(k_n g_n)
=S(t_C)^*j_n(T)=j_n(C).\eqno{\sq}$$

\proc{Remark}\rm In [\Rackspace] we shall relate James bundles to multiple
points of immersions: If $f\co M\to
C$ is a transverse map, then the Poincar\'e duals of the pull-backs 
of the James bundles comprise a self-transverse immersion and its
multiple manifolds.   Using this result, the connection we have given 
of James bundles with classical
James--Hopf invariants overlaps with the results of 
Koschorke and Sanderson [\KS] which also relate 
generalised James--Hopf invariants to multiple points of immersions.
The results overlap  when considering  $\Omega S(S^1)$. Here $\Omega S(S^1)$
is a special case of $\mo C\mo$. In [\KS], $\Omega S(S^1)$ is a special 
case of $\Omega^nS^n(T(\xi))\>,n\ge1,$ where $T(\xi)$ is a Thom space.

\section{James classes and characteristic classes}

In this section we examine the algebra generated by the James classes
and recover (stably) the results of Baues [\Baues].  We also consider
the reduction in ordinary cohomology, which define natural
characteristic classes for a $\sq$--set whose mod 2 reductions are
pulled back from both the Stiefel--Whitney and the Wu classes of $BO$.

The James--Hopf invariants $j_n(C)$ of a $\sq$-set $C$ are classes in
unstable cohomotopy, which stabilise to classes $\gamma_n(C)\in {\Bbb
S}^n(C)$ in stable cohomotopy.  There is a cup product defined for the
unstable classes. Baues [\Baues; page 79] gives a formula for these
cup products for the case $C=T$ (ie for $S^1_\infty$). See Dreckmann
[\Dre; page 42] for a correction.  Here we give a simple geometric
proof of the formula after one suspension.

The (unstable) cup product
$$[SX,S^{m+1}]\times[SX,S^{n+1}]\to [SX,S^{n+m+1}]$$ may be described
in terms of mock bundles as follows.  Let $\xi^m$ and $\eta^n$ be
framed mock bundles embedded in $(X-\{*\})\times \re$. Assume the
projections, $p_\xi$ and $p_\eta$ are transverse.  Define
$E_{\xi\cupprod\eta}=p^*_\xi E_\eta$. Since $E_\eta$ is
$(n+1)$--framed in $X\times \re$, $E_{\xi\cupprod\eta}$ is
$(n+1)$--framed in $E_\xi\times\re$. But $E_\xi$ is $(m+1)$--framed in
$X\times\re$, so using the last of these framing directions to embed
$E_\xi\times\re$ in $X\times\re$ we have an embedding of
$E_{\xi\cupprod\eta}$ which is $(n+m+1)$--framed in $X\times\re$. If
the roles of $\xi$ and $\eta$ are reversed we get a homeomorphic
result with a possible framing change.

However in the case of the James bundles the last framing direction
is parallel to $\re$ in $X\times\re$ and we can see that the two
cup products agree up to permutation of components of blocks and
a sign change. They also agree (after one suspension and up to
sign) with the external
cup product (take the Cartesian product, which is naturally 
$(n+m+2)$--framed in $X^2\times\re^2$ and restrict to the diagonal).
This is seen by rotating the last framing vector through a right
angle into the last $\re$--coordinate.

Now let
$H$ be an $n$ element subset of $\{1,2,\ldots,m+n\}$. Let $s(H)$
be the sign of the shuffle permutation of 
$\{1,2,\ldots,m+n\}$ which moves $H$
to the front and preserves order in both  $H$ and its complement. 
Let $\phi_{m,n}=\sum_H s(H)$ where the sum runs over all n element subsets
$H$ of $\{1,2,\ldots,m+n\}$.

\proc{Proposition}For any $\sq$-set $C$,\quad
$Sj_m(C)\cupprod Sj_n(C)=\phi_{m,n} Sj_{m+n}(C)$.
\endproc
\prf  We prove the result in the case $C=T$, the trivial complex.
The formula then pulls back over
  the constant map
  $t_C\co C\to T$ since the cup product is natural for cubical maps.
To compute the cup product we must make the mock bundle projections
 $p_m\co |J^m(T)|\to|T|$ and
$p_n\co |J^n(T)|\to|T|$ transverse.
Let $b_1=c_{(i_1,\ldots,i_m)}$ and $b_2=c_{(l_1,\ldots,l_n)}$ be elements of
$J^m(T)$ and  $J^n(T)$ respectively  where $c\in T_k$.
 In order to compute the transverse
   intersection of  $p_m |b_1|$ with $p_n |b_2|$ in $|c|$ consider the corresponding sections
    $s_{(i_1,\ldots, i_m)}\co I^{k-m}\to I^k$ and
     $s_{(l_1,\ldots, l_n)}\co I^{k-n}\to I^k$.
The contribution to $j_m(T)\cupprod j_n(T)$ over $|c|$ is given as follows.
 There are two possibilities.
 If $\{i_1,\ldots,i_m\}\cap\{l_1,\ldots,l_n\}=\emptyset$
  then the sections are transverse and they meet in the image
   of $s_{(p_1,\ldots,p_{n+m})}$ where
$\{p_1,\ldots,p_{n+m}\}=\{i_1,\ldots,i_m\}\cup\{l_1,\ldots,l_n\}$.
 Since, strictly speaking, $p$ is a function $\{1,\ldots,n+m\}\to{\Bbb N}$
  we can
define $H=p^{-1}\{i_1,\ldots,i_m\}$. We have determined
 an element $c_{(p_1,\ldots,p_{n+m})}\in J^{n+m}(T)$ 
whose realisation is a  $k-(n+m)$ cube in the total space
 of the mock bundle $\zeta^{n+m}(T)$. But the cup product
  framing differs from the natural framing by the sign of the shuffle $s(H)$. 
This covers the transverse case. \nl
If say $\{i_1,\ldots,i_m\}\cap\{l_1,\ldots,l_n\}=\{q_1,\ldots,q_s\}$ with
$q_1<\ldots<q_s$ then we isotope $p_n$ on $|c_{(j_1,\ldots,j_n)}|$
 by increasing the
 $q_1^{\rm th}$ coordinate in $I^k$ from $\ha$ to say ${\scriptstyle {3\over4}}$,
  then the intersection is empty and there is no contribution
   to $j^m\cupprod j^n$ in this case. 

At this point we have proved the result provided we can cancel off
contributions of opposite sign and reorder the pieces.  For this
we need to use the one extra dimension where there are obvious cobordisms
which achieve this result. \qed

\proc{Lemma}For $m,n\in{\Bbb N}$
$$\phi_{m,n}=\cases{{[(m+n)/2]\choose [n/2]} &if $nm$ even,\cr
0,&otherwise\cr}$$
In particular for all $m,n\quad 
\phi_{2m,2n}={m+n\choose m}$,
$\phi_{m,n}={m+n\choose m} \ \hbox{mod}\  2$,
$\phi_{1,1}=0, \phi_{1,2n}=1$.
\endproc
\prf Clearly $\phi_{0,n}=\phi_{m,0}=1$. Set $\phi_{m,n}=0$
 if either $m$ or $n$
is negative. For $\phi_{m,n}$ the sum can be split up,
according as the least element of $H$ is or is not $1$,
 giving the recurrence relation:
$$\phi_{m,n}=\phi_{m,n-1}+(-1)^n\phi_{m-1,n}.$$
The  result now follows by induction.\qed

The computation of the algebra generated by the stabilisation
of the James--Hopf invariants easily follows. We shall carry out
this computation after introducing the characteristic classes of
a $\sq$--set, which are the corresponding
cohomology classes.

\sh{The characteristic classes of a $\tsq$--set}

Let $C$ be any $\sq$--set.
The $n$--{\sl th characteristic class} of $C,\quad V_n(C)\in H^n(C)$ 
is the class of the unit cocycle: value $1$ on each $n$--cube.\endproc

Let ${\Bbb S}$ denote the sphere spectrum and let 
$\ga_n(C)\in{\Bbb S}^n(C)$ denote the stabilisation of $j_n(C)$.
The element $\ga_n(C)$ is represented by $E_{\zeta^n}$ framed
in $\mo C\mo\times\re^{\infty}$. Let $h\co {\Bbb S}^*\to H^*$
be the Hurewicz transformation.

\proc{Proposition} For any $\sq$--set $C$\quad $h(\ga_n(C))=V_n(C)$.\endproc

\prf At the mock bundle level $h$ can be described 
simply by counting, with sign, the components
of the $0$--dimensional blocks. The result follows
since the $0$--dimensional blocks of  $\ga_n$ consist
of  a single point over each $n$-cube with standard framing.\qed

Denote the exterior algebra, over $\Z$, on one generator $x$ by
$E(x)$, then $x^2=0$ and elements of $E(x)$ can be written as
$p+q.x,\quad p,q\in \Z$. Denote the divided polynomial algebra on one
generator $y$ by $D(y)$, in other words there are elements
$y_n,\>n\ge0,$ with $y_0=1$ and $y_1=y$ such that $y_n.y_m={n+m\choose
m}y_{n+m}$.  An element of $D(y)$ is then just a linear combination of
the $y_n$'s.

\proc{Corollary}The elements $\ga_n(T),\> n\geq 0,$ generate a subalgebra of 
${\Bbb S}^*(|T|)$ of form $E(\ga_1(T))\otimes D(\ga_2(T))$.
 In particular any element of the subalgebra can be
 uniquely written as a linear combination of the elements $\{1,\ga_1\ga_{2n},\ga_{2n}\}_{n\geq0}$.\endproc

\prf  The result follows from the proposition and 
the lemma and the fact noted above that the Hurewicz map 
$h\co {\Bbb S}^*(|T|)\to H^*(|T|)$ is non trivial on the
 $\ga_n, \>n\geq0$. In fact it
 follows that  $h$ restricted to the subalgebra is an isomorphism.\qed

\sh{The $\Z_2$--characteristic classes}

Let $v_n$ denote the mod 2 reduction of our characteristic class $V_n$. 
We can identify $v_n(T)$ as the Stiefel--Whitney class of a vector bundle. Background for what follows  can be found in [\SMah] and [\Mah]. Let $BO$ be the classifying space for the infinite orthogonal group, and let $B^2O$ be
its delooping; so $\Omega B^2O\simeq BO$. Let $\omega\co S^2\to B^2O$ be a generator of $\pi_2B^2O\cong\Z_2$.
By taking $\Omega\omega$ we get a map $s\co|T|\to BO$ with the property that $s$ restricted to the $1$-skeleton classifies the M\"obius band.

\proc{Proposition} $w_n(s^*\gamma)=w_n(s^*\gamma^{-1})=v_n(T),\>n\geq0$,
where $\gamma$ is the universal (virtual)
bundle.\endproc
\prf    Consider the $n$--fold product $S_n=S^1\times\ldots\times S^1$. Let
$m\co S_n\to |T|$ be given by multiplication. Then 
$m^*\co H^n(|T|)\to H^n(S_n)$ is an isomorphism. Since $s$ is an $H$-map
$w(m^*s^*\gamma)=\prod_{i=1}^n(1+a_i)$ where each
$a_i$ is the generator of $H^1(S^1,\Z_2)$, and so $w_n(s^*\gamma)\not=0$.
Since $m^*s^*\gamma$ is a product of (stable) line bundles over
 $S^1$'s $m^*s^*\gamma=m^*s^*\gamma^{-1}$. \qed

Let $\phi_C\co |C|\to BO$ be the composition $s\circ t$ where 
$t\co|C|\to |T|$ and let $\ga_C=\phi_C^*\ga$.
We call $\ga_C$ the {\sl characteristic bundle} over $C$.

\proc{Corollary}For any $\sq$--set $C$, 
\ $w_n(\ga_C)=v_n(C),\>n\geq0$.
\qed

\rk{Remark}
Recall that the (total) Stiefel--Whitney class, $w(\xi)$, of a vector bundle
$\xi$ is defined by $w(\xi)\cupprod U=Sq\>U$ where $U$ is the Thom class.
Similarly the Wu class, ${\rm wu}(\xi)$, is given by ${\rm wu}(\xi)\cupprod U=\chi
Sq\>U$. From this and the Cartan formula we get the well known
relation $Sq\>{\rm wu}(\xi)=w(\xi^{-1}).$  On $S_n$ we have $Sq=Id$ and $m^*s^*\ga=m^*s^*\ga^{-1}$
so both in the proposition and the corollary
Whitney and Wu classes coincide.

We shall continue to study the characteristic bundle and the
associated map\break $\phi_C\co |C|\to BO$ in the next section, 
where we use it
to define the concept of a $C$--oriented manifold and to define an
associated generalised cohomology theory.

\section{Homology Theories from $\square$-sets}

In this section we describe how a $\sq$--set may be used to 
construct homology and cohomology theories.
More precisely we shall define a spectrum $MC$ for any $\sq$--set $C$ 
which then defines a generalised homology theory $MC_*$.  This theory
has a good description as a bordism theory in terms of manifolds and
labelled diagrams.  The description is particularly nice when $C$ is
the classifying space of a rack.  These theories have natural
operations, which as we shall see are analogous to ${\cal X}Sq$ 
operations in ordinary cohomology.

In the case $C=T$, the trivial $\sq$--set, the spectrum $MT$ 
is a ring spectrum and
is one of a family studied by Mahowald [\Mah].

\sh{The spectrum $MC$}

We recall from section 4 that there is a map $s\co\mo T\mo\to BO$
(where $T$ is the trivial $\sq$--set with just one cell in each
dimension and $BO$ is the classifying space for stable vector bundles) and
this defines a map $\phi_C=s\circ t_C\co \mo C\mo\to BO$ for any
$\sq$--set $C$, where $t_c\co C\to T$ is the canonical $\sq$--map.  
The map $\phi_C$ has the property that it pulls
both the universal Stiefel--Whitney and Wu classes back to the
mod 2 reductions of the characterstic classes of $C$, see 4.6.

Let $\ga$ be the universal (virtual) bundle over $BO$ and $\ga_C$
(the characteristic bundle) 
its pull-back over $C$ by $\phi_C$.  The spectrum $MC$ is defined to
be the associated Thom spectrum to $\ga_C$.  It defines a generalised
homology theory, which we call {\sl $C$--bordism}, and
denote $MC_*$ and a corresponding theory of {\sl $C$--cobordism}
denoted $MC^*$.

\sh{Connection with standard spectra}

We next see how the new spectrum fits with the sphere spectrum (which
classifies stable cohomotopy) and the Eilenberg-MacLane spectrum
which classifies ordinary (mod 2) cohomology.

The universal Thom class determines a map of spectra $MO\to H\Z
_2$, where $H\Z_2$ is the mod 2 Eilenberg-MacLane spectrum. Hence we
get  a natural map $\theta_C\co MC\to H\Z_2$.
Let ${\Bbb S}$ denote the sphere spectrum, and let ${\Bbb S}\Z_2$ denote
the sphere spectrum with mod 2 coefficients. 
\proc{Proposition}For any $\sq$--set $C$:\items
\item{\rm(1)}The Hurewicz map factors as ${\Bbb S}\to MC\to H\Z_2$,
\item{\rm(2)}If there is an element $x\in C_1$ with $\d_1^0x=\d_1^1x$
then the mod 2 Hurewicz map factors as ${\Bbb S}\Z_2\to MC\to H\Z_2$.
\enditems\endproc

\rk{Note}The condition in part (2) is always satisfied for rack spaces.

\prf   The sphere spectrum ${\Bbb S}$ may be regarded as the Thom spectrum
 pulled back from $BO$ by the inclusion of a point. The first result follows
from naturality of $\theta$ by choosing a vertex in $\mo C\mo$.

For the second part let $D$ be the $\sq$--set which is the
circle with the obvious $\sq$--set structure, ie 
$D$ has a single vertex, a single edge and no higher cells.
From the property that $\phi_D$ pulls back the Stiefel--Whitney
class to the mod 2 characteristic class of $D$, it can be seen that
$\phi_D\co \mo D\mo\to BO$ classifies the (stable) M\"obius band
over the circle $\mo D\mo$. The Thom complex of the M\"obius band is the 
projective plane and $MD$ is its suspension spectrum, ie ${\Bbb S}\Z_2$.
The transformation $\theta_D\co MD\to H\Z_2$ can now be identified with
the inclusion of the (stable) $1$--skeleton ie the mod 2 Hurewicz map.
The result now follows from naturality of $\theta$ by considering 
the inclusion of $D$ in $C$ given by mapping the edge of $D$
to $x$.\qed

\sh{Geometric description of the cycles of $MC$}

The generalised homology theory defined by $MC$ is the bordism
theory defined by a certain class of manifolds.  We shall describe
this class.

We need to define what it means for a manifold to have its normal
bundle {\sl twisted} by a self-transverse immersed submanifold
of codimension 1.  For simplicity we describe the
case of an embedding first.  So let $Q\subset M$
be a codimension 1 submanifold and assume for simplicity that
$Q$ is equipped with a bicollar
neighbourhood $N(Q)$.  Denote by $\rho/I$ the {\sl M\"obius bundle}
which is the 1--dimensional vector bundle constructed by
gluing the trivial bundles over $[0,\ha]$ and $[\ha,1]$ by 
multiplication by $-1$ in the fibre over $\ha$.  Thus although 
$\rho/I$ is the trivial
bundle it is equipped with opposite trivialisations on the two
halves of $I$.  Denote  by $\hat\rho/I$ the stable version.
We call $\hat\rho/I$ the (stable) {\sl M\"obius twist}.

Now construct the stable bundle $\hat\rho/M$ by taking the
trivial bundle over $Q^c=\bar{M-N(Q)}$ and gluing in $\hat\rho/I$
on each collar line.  Alternatively we can take the trivial
line bundle on $M-Q$ and glue at $Q$ by multiplication by $-1$
and then stabilise.
We say that the stable normal bundle $\nu/M$ of $M$ is {\sl twisted} by
$Q$ if it is identified with $\hat\rho/M$.

The concept of twisting by a self-transverse immersion is similar 
except that at $k$--tuple
points there are $k$ mutually perpendicular M\"obius twists taking place
corresponding to the directions perpendicular to the $k$ sheets.

For the detailed definition we shall contruct a stable M\"obius
twist bundle over $I^n$ for each $n$.  To do this we need a
function $v\co J^1(I^n)\to S^{n-1}$ with an orthogonality
property.  Denote $v\mo I^n_{(i)}$ by $v_i$ which we think of as a 
variable unit vector in $\re^n$.  We require that, for each point 
$p\in I^n$ covered by points
$p_{i_1},\ldots p_{i_t}$ in $I^n_{(i_1)},\ldots I^n_{(i_t)}$,
the corresponding vectors $v_{i_1},\ldots v_{i_t}$ should
form an orthonormal $t$--frame.  We also require that $v$ commutes
with face maps, ie the diagram commutes
$$
\matrix{
J^1(I^{n-1})&{\buildrel v \over\a}&S^{n-2}\cr
&&\cr
\downarrow\!{\scriptstyle  J^1(\delta_k^\ep)}&&\downarrow\!{\scriptstyle i}\cr
&&\cr
J^1(I^{n})&{\buildrel v\over \a}&S^{n-1}\cr}
$$
where $i\co S^{n-2}\to S^{n-1}$ is the usual inclusion.
\def\int{\hbox{\rm int}}

We define $v$ by induction on $n$ and suppose that 
$v\co J^1(I^{n-1})\to S^{n-2}$ is already defined.  Then the
commuting diagram defines $v$ on $\d I^n$.  We extend to the
interior of $I^n$ by using the cubical subdivision of $I^n$
determined by the image of $J^1(I^n)$ in $I^n$, which cuts
$I^n$ into $2^n$ subcubes.  Note that $n-i$ sheets of the
image of $J^1$ meet at each $i$--cube in $\int I^n$ in this 
subdivision (see the diagram near the end of section 1).  
We define $v$ by induction over skeleta.
At the unique 0--cell in $\int I^n$ we choose $v_1,\ldots,v_n$
to be any orthonormal $n$--frame.  In general we have an orthonormal
$(n-i)$--frame defined over the boundary of a typical $i$--cell
and we extend over the $i$--cell using the fact that 
$\pi_{i-1}(V_{n,n-i})=0$.

This completes the definition of $v$.  Notice that the choices made
in the definition are homotopic if $S^{n-1}$ is replaced by $S^n$ 
and hence the stable
twisted bundles constructed below do not depend on these choices.

We now construct the $n$--dimensional vector bundle $\rho/I^n$
by taking the trivial bundle off the image of $J^1$ and gluing
by mapping $v_i\mapsto -v_i$ as the $i^{\rm th}$ sheet (the image of
$I^n_{(i)}$) is crossed and leaving the orthogonal complement
of $v_i$ fixed.  Near a point $p\in I^n$ covered by points
$p_{i_1},\ldots p_{i_t}$ in $I^n_{(i_1)},\ldots I^n_{(i_t)}$
the various trivial bundles are being glued by 
reversing the signs of the appropriate subset of
$v_{i_1},\ldots v_{i_t}$.  By orthogonality, these $t$ simultaneous
gluings are all independent.  We denote the stable version of
$\rho/I^n$ by $\hat\rho^n$.  This is the {\sl $n$--fold M\"obius twist}
over $I^n$.

By construction, M\"obius twists have the property that for 
each face map $\la\co I^t\to I^p, \la^*(\hat\rho^p)=\hat\rho^t$.
This implies that if $X$ is a $\sq$--space then we can construct
a bundle $\hat\rho/\mo X\mo$ by gluing together the bundles $X_n\times
\hat\rho^n$ mimicking the definition of $\mo X\mo$.  Note that
$\hat\rho/\mo X\mo$ is trivialised over $X_0$ and at $X_n$ there
is an $n$--fold M\"obius twist taking place.

\sh{Neighbourhood systems and transversality}

Now suppose that $Q$ is the image of a self-transverse [\LS] immersed
submanifold of $M^m$ of codimension 1.  We can construct a {\it
neigbourhood system} $N(Q)$ as follows.  We choose a cube
neighbourhood of each $m$--tuple point, which (using
self-transversality) can be assumed to meet nearby sheets in the $m$
central $(m-1)$--subcubes (ie the images of $J^1$).  Then, working in
the complement of the interiors of these cubes we choose $I^{m-1}$
bundle neighbourhoods of the 1--dimensional $(m-1)$--tuple points
which meet the cubes near the $m$--tuple points in faces.  We continue
to choose inductively $I^q$ bundle neighbourhoods of the $q$--tuple
points meeting nearby sheets in central subcubes and higher multiple
points in faces.  We finish by choosing (trivial) $I^0$ bundle
structures at the $0$--tuple points (ie the points of $M-Q$).  Thus
the final neighbourhood system decomposes $M$ into $I^n$ bundles for
varying $n$, glued together along the analogue of faces.

We need the following theorem.

\proc{Transversality theorem}Let $C$ be a $\sq$--set and $f\co M^m
\to\mo C\mo$ a map.  Then $f$ is homotopic to a map $f'$
(called {\sl transverse}) which has the following property.  There is
a self-transverse immersed submanifold $Q\subset M$ of codimension 1
and a neighbourhood system $N(Q)$ which is in fact the realisation of
a $\sq$--space such that $f'\co M\to \mo C\mo$ is the realisation of a
$\sq$--map.\rm

\prf  The proof is similar to the proof of transversality for CW
complexes given in [\BRS; Chapter 7].  We construct $f'$ by a sequence
of changes using a downwards induction on the dimension of skeleta of
$C$ containing the image of $f$.  By the usual abuse of notation we
will continue to call the amended maps $f$.  By cellular approximation
we can assume that the image of $f$ lies in the $m$--skeleton of $C$.
Use smooth transversality to make $f$ transverse to the centres of the
$m$--cubes of $C$ and then by expanding a small neighbourhood of each
centre onto the whole cube, we can assume that the preimage of each
$m$--cube is a collection of $m$--cubes in $M$ and that $f$ maps each
by a standard indentification to the corresponding $m$--cell of $C$.
Remove the interiors of these $m$--cubes and work in their complement.
The map $f$ restricted to the boundary maps to the $(m-1)$--skeleton
of $C$ and is transverse to centres of $(m-1)$--cells.  Consider one
such centre point $p$ in an $(m-1)$--cell $R$.  By relative
transversality we can assume that the preimage of $p$ is a 1--manifold
with normal bundle mapped to a neighbourhood of $p$ in $R$.  By
expanding a small neighbourhood of $p$ onto the $(m-1)$--cell we can
assume that the normal bundle is fibred by copies of $R$ mapped by the
identity to $R$.  We now excise the interiors of these bundles and we
are mapping to the $(m-2)$--skeleton and we repeat the construction
using the centres of $(m-2)$--cells.  Continue using the centres of
$(m-p)$--cells, $p=3,\ldots,m$. The end result of this process is the
required transverse map.  Let $Q$ be the self-transverse immersed
submanifold which comprises all the codimension 1 central subcubes in
all the fibres (the cubes of $J^1(R)$ as $R$ varies).  For example
these would be the planes with bold borders in figure \JaComp).  The
cube bundles of varying dimensions that have been constructed form the
required neighbourhood system of $Q$.
\qed

\rk{Remark}We can describe $\hat\rho/I^n$ as an explicit subbundle
of the infinite trivial bundle in a similar way to the analogous
construction in [\SMah].  To be precise we can embed
$\hat\rho^n$ as a subbundle of $I^n\times\re^n\times\re^n$ as follows.
Choose a small neighbourhood system for the image of $J^1(I^n)$
in $I^n$, regarded as usual as made of (trivial) cube bundles
of varying dimensions with fibres which we shall call ``small''
cubes (to distingush them from the bigger cubes which comprise
$I^n$ and $J^1(I^n))$.   Without loss assume that the function
$v$ is constant on the preimage of each small cube.  Think
of the neighbourhood system as covered by double
collars on the $I^n_{(i)}$ where each collar line lies
in a small cube parallel to an edge.  Now over each such
line (corresponding to $p_i\in I^n_{(i)}$ say) turn the vector
$v_i$ over in the plane determined by $v_i$ and $\bar v_i$ (which
is the copy of $v_i$ in the other copy of $\re^n$) in other words make
an explicit M\"obius bundle in this plane over the line.  Where vectors in
two or more such bundles lie over the same point of $I^n$
they are contained in perpendicular planes, so there is no
interference.  The construction of $\hat\rho/\mo X\mo$, where
$X$ is a $\sq$--space, now yields $\hat\rho$ as an explicit subbundle
of the infinite trivial bundle.

\sh{$C$--manifolds}

Now suppose we are given a transverse map $f\co M\to C$
where $C$ is a $\sq$--set.  Thus we have a self-transverse
immersed submanifold $Q$ in $M$ and an identification of
$M$ with the realisation of a $\sq$--space structure $X$ 
on a neighbourhood system of $Q$.  We say that the stable 
normal bundle $\nu/M$ of $M$ is {\sl twisted} 
by $Q$ if it is identified with $\hat\rho/\mo X\mo$.

We can now define a $C$--manifold.  
This is a manifold $M$ and a transverse map $f\co M\to C$,
such that the stable normal bundle $\nu/M$ of $M$ is twisted 
by the corresponding self-transverse immersed submanifold. 
There is an obvious concept of {\sl $C$--cobordism}
between $C$--manifolds and then we can define the $C$--bordism groups 
of a space $X$ by mapping $C$--manifolds and cobordisms into $X$ in the
usual way.  There is a dual concept of $C$--cobordism given by
considering mock bundles with fibres $C$--manifolds.

\proc{Theorem}The theory of $C$--bordism given by $C$--manifolds
coincides with the theory $MC_*$ defined by $MC$.  Similary the
two cobordism theories coincide.\rm

Before proving the theorem we shall recall another theory
considered by Mahowald [\Mah] for which a geometric
description is given in [\SMah].

\sh{The Mahowald spectrum}

The spectrum $MA$ (another of the spectra defined by Mahowald
[\Mah]) is pulled back from $MO$ by $\Omega^2\nu$ where $\nu\co
S^3\to B^3O$ is a generator.  Mahowald [\Mah] and Priddy [\Priddy]
have proved that the composition $MA\to MO\to H\Z_2$ is an
equivalence.   

\proc{Proposition}$MC\to MO$ factors as $MC\to MA\to MO$.\rm

We shall prove the proposition and the theorem together.
We shall need to recall the geometric
description of $MA_*(-)$ given in [\SMah].  We shall translate
the language of [\SMah] into the language of this paper.
Roughly speaking, ``good position'' as in [\SMah] corresponds
to the neighbourhood systems used here.

Consider a framed self-transverse immersion of codimension 1
$Q\to M$.  A neighbourhood system of $Q$ can be constructed as
described above.  Although this does not necessarily give $M$ a
$\sq$--space structure, twisting by $Q$ can be
described in a very similar way.  The details
can be obtained from [\SMah].

An {\sl $A$--manifold} is a manifold $M$ with a framed self-transverse 
immersion of codimension 1 $f\co Q\to M$ covered by an embedding
$\tilde f\co Q\to M\times\re^2$ and such that $\nu_M$ is twisted
by $Q$.   The following result is proved in [\SMah].

\proc{Theorem}The bordism theory determined by $A$--manifolds is the
theory $MA_*$.  Similarly the cohomology theories coincide. \qed

\sh{Proofs of theorem 5.3 and proposition 5.4}

An element of the homology group $MC_m(Y)$ is
determined by a closed $m$--manifold $M^m$ together with a 
map  $f\co M\to Y$ and an isomorphism  $F\co\nu_M\to\phi_C^*\ga$
where $\nu_M$ is the stable normal bundle of $M$ and $\ga$
is the universal bundle over $BO$.  Now $F$ restricts to a map 
$M\to \mo C\mo$ which may be assumed to be transverse
by 5.2 and determines a codimension 1 self-transverse immersed
submanifold $Q$ say.

It remains to prove that the isomorphism $F\co\nu_M\to\phi_C^*\ga$
is equivalent to $\nu_M$ being twisted by the diagram.  To see
this we shall use the spectrum $MA$ and the geometric description
given above.

There is a well known
configuration space model $C_2(S^1)$ for $\Omega^2S^3$.
$$C_2(S^1)=(\coprod C_{2,k}\times_{\Sigma_k}I^k)/\sim$$
where $C_{2,k}$ denotes ordered subsets of $\re^2$ containing $k$
points and
$$[x_1,\ldots, x_k,t_1,\ldots,t_k]\sim[x_1,\ldots,\hat x_i,\ldots,
x_k,t_1,\ldots,\hat t_i,\ldots, t_k]\hbox{\ \ if\ \ }t_i\in\{0,1\}.$$

The space of unordered sets $C_{2,k}/\Sigma_k$ can be identified with
the classifying space of the braid group $Bbr_k$. There is a
filtration of $C_2(S^1)$:
$$*\subset F_1C_2(S^1)\subset\ldots\subset F_kC_2(S^1)\ldots
C_2(S^1)\simeq \Omega^2S^3.$$
where $F_k$ is obtained from $F_{k-1}$ by attaching the bundle of
cubes $C_{2,k}\times _{\Sigma_k}I^k\to Bbr_k$ along its sphere bundle 
$C_{2,k}\times _{\Sigma_k}\d(I^k)\to Bbr_k$. 

There is a similar model for $\Omega S^2$. Simply replace
$C_{2,k}$ by $C_{1,k}$. There are maps 
$\phi_k\co I^k\to C_{1,k}\times _{\Sigma_k}I^k$ given by
$(t_1,\ldots,t_k)\mapsto [(0,0),\ldots,(k,0),t_1,\ldots,t_k]$.
These maps combine to give a homotopy equivalence
$|T|\to C_1(S^1)$.  The composition
$$|T|\to C_1(S^1)\subset C_2(S^1)\simeq 
\Omega^2S^3{\buildrel\Omega^2\nu\over \longrightarrow}BO$$
is our standard map $s\co |T|\to BO$, and so we have
$$MC\to MA\to MO\to H\Z_2$$
and it follows that the map $M\to BO$ factors via $MA$.

The required result now follows from theorem 5.5 above after
noting that the descriptions of twisting given above and
in [\SMah] are effectively the same. \qed

\rk{Remark}
The maps ${\Bbb S}\to MC$ and $MC\to MA\simeq H\Z_2$ can now be described
geometrically. To get from ${\Bbb S}_*(-)$ to $MC_*(-)$ simply take
the map to $C$ to be the constant map to the chosen vertex.  To
pass from $MC_*(-)$ to $MA_*(-)$ forget the map to $C$, remember only
the codimension 1 self-transverse submanifold and allow
normal bundles to be classified by the cube bundles over $Bbr_k$.

\sh{Restrictions on normal bundles}

The manifolds which can occur in the definition of $MC_n(X)$ are restricted
by the condition on the normal bundles. In particular we have the
following result.

\proc{Proposition} If $M$ is a $C$--manifold
then ${\rm wu}(M)=w(M)$, and $w_i(M)=0$ for $2i>n$.\endproc

\prf  Recall that ${\rm wu}(M)={\rm wu}(\nu(M))$ and $w(M)=w(\tau(M))$. The
result now  follows from the remark after 4.6 and the fact that
Wu classes of manifolds vanish above half dimension. \qed

\sh{Rack and link homology theories}

In the case when $C=BX$ where $X$ is a rack, the geometric
description of $MC_*$ can be simplified.  A transverse map
$f\co M\to BX$ is equivalent to a framed codimension 2
submanifold of $M\times\re$ together with a representation of 
the fundamental rack in $X$ (see [\Class]).  If $X$ is the
fundamental rack of a link $L$, then we can regard the resulting
theory as $L$--bordism, and this provides many more new invariants 
of links by calculating the values of this theory on test spaces.

There is also a ``truncated theory'' given by using only
the 2--skeleton of the rack space and this also makes good
sense for any trunk as well.
The resulting spectrum can be identified with ${\Bbb S}\wedge X$ 
where $X$ is the Thom complex of a plane bundle over the $2$-skeleton.

\sh{Operations}

The geometric description of $C$--manifolds given above allows us 
to define operations on these theories, which cover
inverse Steenrod square operations in $\Z_2$--cohomology.  

\proc{Theorem}Let $C$ be a $\sq$--set.  There is
an operation $MC_n\to MD_{n-k}$ where $D=J^k(C)$ and a similar
operation $MC^n\to MD^{n+k}$.  These operations correspond to 
a map of spectra $MC\to \Sigma^kMJ^k(C)$. 

In case $C=T$, the trivial
$\sq$--set,  we can compose with $J^k(T)\to T$ to get a map 
$MT\to \Sigma^kMT$ which fits into the following commutative diagram
$$\matrix{
MT & \a & H\Z_2\cr
       &     &           \cr
\downarrow &  &\ \ \ \ \ \downarrow {\cal X}Sq^k\cr
&     &           \cr
\Sigma^kMT & \a & \Sigma^kH\Z_2\cr
}$$
where horizontal maps are the map $MT\to MA \cong H\Z_2$ considered
above and its $k$--fold suspension.

\prf  We need a description of the $k$--tuple points of the
immersion $Q$ determined by a transverse map $f\co M\to \mo C\mo$
given in [\Class; section 2], namely that there is an induced
transverse map $f^{(k)}\co Q^{(k)}\to J^k(C)$ where $Q^{(k)}$
is the manifold covering the $k$--tuple points of $Q$.
It can be seen that the condition that the normal bundle is twisted
by the immersion is inherited by the induced immersion in $Q^{(k)}$
and this defines the operations. 

Now the map $MT\to MA \cong H\Z_2$ is induced by the universal Thom
class of $MT$ and the commutativity of the diagram then follows from
the results of [\SMah], see in particular remark 3.6. \qed

\rk{Remark} The operation $MT\to \Sigma^kMT$ defined by theorem 6.6 
coincides with one of the operations defined by Mahowald on p552 of
[\Mah].

\references

\bye

%% file: newinsert.tex
\chardef\newinsCatAt\the\catcode `\@
\catcode `\@=11
%
%%%%%%%%%%%% Corrected insert macros for plain.tex %%%%%%%%%%%%%%%%%
%
%  New skipamounts:
%
\newskip\insertskipamount\newskip\inserthardskipamount
\insertskipamount 12pt plus2pt  %Redefined (CPR) to suit %default 6pt plus 2pt
\inserthardskipamount 4pt       %the G and T style       %default 6pt
\def\insertskip{\vskip\insertskipamount}
%
%  Save and restore \lastskip:
%
\newskip\LastSkip
\def\SaveLastSkip{\LastSkip\lastskip}
\def\RestoreLastSkip{\nobreak\vskip-\LastSkip\vskip\LastSkip}
%
%  Larry Siebenmann's test for split topinserts:
%
\newcount\SplitTest%        will be set to -1 if a topinsert has split
\def\SetSplitTest{\SplitTest\insertpenalties
  \insert\topins{\floatingpenalty1}%
  \advance\SplitTest-\insertpenalties}
%
%  From here on we modify definitions in plain.tex.
%
% Redefine \midinsert to convert to \topinsert if a topinsert has been
% split, to prevent midinserts getting out of order (cf. TeXbook Exercise
% 15.5). As in plain.tex, a \midinsert still converts to a \topinsert
% (which then splits) if the insert is too big for current page.
%   Was:    \def\midinsert{\@midtrue\@ins}
\def\midinsert{\par
 \SaveLastSkip\penalty-150\SetSplitTest\RestoreLastSkip
 \ifnum\SplitTest=-1
  \@midfalse\p@gefalse\else\@midtrue\fi\@ins}
% Redefine \@ins to add \inserthardskipamount of glue above.
%   Was:  \def\@ins{\par\begingroup\setbox\z@\vbox\bgroup}
\def\@ins{\par\begingroup\setbox\z@\vbox\bgroup%
  \vglue\inserthardskipamount}
% Changes to \endinsert of plain.tex 3.0:
% - Use \insertskipamount instead of \bigskipamount throughout.
% - Use larger of previous skip and insertskip before middle insert.
% - Add \nointerlineskip to avoid unwanted extra 1pt skip.
% - Save and restore lastskip when an insert floats.
\def\endinsert{\egroup % finish the \vbox
  \if@mid \dimen@\ht\z@ \advance\dimen@\dp\z@
    \advance\dimen@\insertskipamount%            was 12pt (wn)
    \advance\dimen@\pagetotal\advance\dimen@-\pageshrink
    \ifdim\dimen@>\pagegoal\@midfalse\p@gefalse\fi\fi
% Next 3 lines replace:  \if@mid \bigskip\box\z@\bigbreak (wn)
  \if@mid%
    \ifdim\lastskip<\insertskipamount\removelastskip\insertskip\fi
    \nointerlineskip\box\z@\penalty-200\insertskip
  \else%
    \SaveLastSkip%                                  added (wn)
    \insert\topins{\penalty100 % floating insertion
    \splittopskip\z@skip
    \splitmaxdepth\maxdimen \floatingpenalty\z@
    \ifp@ge \dimen@\dp\z@
    \vbox to\vsize{\unvbox\z@\kern-\dimen@}% depth is zero
    \else \box\z@\nobreak\insertskip\fi}% was \bigskip\fi (wn)
    \RestoreLastSkip%                               added (wn)
   \fi\endgroup}
%%%%%%%%%%%%%%%%%% Done correcting insert macros %%%%%%%%%%%%%%%%%%%
%
\catcode `\@=\newinsCatAt